\documentclass[12pt]{amsart}
\font\emailfont=cmtt10

\usepackage{amsmath,amsthm,amsfonts,amscd,flafter,epsf}

\newcommand\commentable[1]{#1}

\newcommand\Id{\mathrm{Id}}

\newcommand{\rk}{\mathrm{rk}}
\newcommand{\HF}{HF}

\newtheorem{theorem}{Theorem}[section]

\newtheorem{defn}[theorem]{Definition}

\def\endproof{\relax\ifmmode\expandafter\endproofmath\else
  \unskip\nobreak\hfil\penalty50\hskip.75em\hbox{}\nobreak\hfil\bull
  {\parfillskip=0pt \finalhyphendemerits=0 \bigbreak}\fi}
\def\endproofmath$${\eqno\bull$$\bigbreak}
\def\bull{\vbox{\hrule\hbox{\vrule\kern3pt\vbox{\kern6pt}\kern3pt\vrule}\hrule}}

\newcommand{\Q}{\mathbb{Q}}
\newcommand{\R}{\mathbb{R}}

\newcommand{\C}{\mathbb{C}}

\newcommand{\Z}{\mathbb{Z}}

\newcommand{\OneHalf}{\frac{1}{2}}

\newcommand{\Zmod}[1]{\Z/{#1}\Z}

\newcommand{\grad}{\vec\nabla}

\newcommand{\cm}{\cdot}

\newcommand{\Nbd}[1]{{\mathrm{nd}}(#1)}
\newcommand{\nbd}[1]{\Nbd{#1}}

\newcommand{\ModSWfour}{\mathcal{M}}
\newcommand{\ModFlow}{\ModSWfour}

\newcommand{\SpinC}{{\mathrm{Spin}}^c}

\newcommand\abuts\Rightarrow
\newcommand\Sym{\mathrm{Sym}}

\newcommand\spinccanf{k}

\newcommand\HFpRed{\HFp_{\red}}
\newcommand\HFmRed{\HFm_{\red}}

\newcommand\RP[1]{{\mathbb{RP}}^{#1}}

\newcommand\relspinc{\underline{\spinc}}
\newcommand\ThreeCurveComp{\Sigma-\alpha_{1}-\ldots-\alpha_{g}-\beta_{1}-\ldots-\beta_{g}-\gamma_1-...-\gamma_g}

\newcommand\Filt{\mathcal F}
\newcommand\HFinfty{\HFinf}
\newcommand\CFinfty{\CFinf}

\newcommand\x{\mathbf x}
\newcommand\w{\mathbf w}

\newcommand\y{\mathbf y}

\newcommand\ModSphere{\ModFlow\left({\mathbb S}\longrightarrow 
\Sym^{g-1}(\Sigma_{1})\times \Sym^2(\Sigma_{2})\right)}
\newcommand\ModSpheres\ModSphere
\newcommand\CF{CF}

\newcommand\CFa{\widehat{CF}}
\newcommand\CFp{\CFb}
\newcommand\CFm{\CF^-}

\newcommand\HFpred{\HFp_{\rm red}}

\newcommand\HFmred{\HFm_{\rm red}}
\newcommand\HFred{\HF_{\rm red}}

\newcommand{\red}{\mathrm{red}}

\newcommand\HFp{\HFb}

\newcommand\HFm{\HF^-}
\newcommand\CFinf{CF^\infty}
\newcommand\HFinf{HF^\infty}
\newcommand\CFb{CF^+}
\newcommand\HFa{\widehat{HF}}
\newcommand\HFb{HF^+}
\newcommand\gr{\mathrm{gr}}
\newcommand\Mas{\mu}
\newcommand\UnparModSp{\widehat \ModSp}
\newcommand\UnparModFlow\UnparModSp
\newcommand\Mod\ModSp

\newcommand\spin{\mathfrak s}

\newcommand{\spinc}{\mathfrak s}

\newcommand{\spinct}{\mathfrak t}

\newcommand\Real{\mathrm Re}

\newcommand\ModMaps{\mathcal M}
\newcommand\ModSp\ModMaps

\newcommand\Ta{{\mathbb T}_{\alpha}}
\newcommand\Tb{{\mathbb T}_{\beta}}
\newcommand\Tc{{\mathbb T}_{\gamma}}

\newcommand\Strip{\mathbb{D}}

\newcommand\alphas{\mbox{\boldmath$\alpha$}}

\newcommand\betas{\mbox{\boldmath$\beta$}}
\newcommand\gammas{\mbox{\boldmath$\gamma$}}

\newcommand\Fm[1]{F^{-}_{#1}}

\newcommand\Fp[1]{F^{+}_{#1}}

\newcommand\Finf[1]{F^{\infty}_{#1}}

\newcommand\SpinCCobord{\theta^c}

\newcommand\Fc{F^\circ}

\newcommand\HFc{\HF^\circ}

\newcommand\Dual{\mathcal D}
\newcommand\Duality\Dual

\newcommand\HMfrom{\widehat{\mathrm{HM}}}
\newcommand\HMred{\overline{\mathrm{HM}}}
\newcommand\HMto{\check{\mathrm{HM}}}

\newcommand\InjMod{\mathcal T}

\newcommand\spincrel\relspinc

\newcommand\Vertices{\mathrm{Vert}}

\newcommand\HFK{HFK}

\newcommand\HFKa{\widehat\HFK}

\commentable{

\title[{Heegaard diagrams and holomorphic disks}]
{Heegaard diagrams and holomorphic disks}

\author[Peter Ozsv{\'a}th]{Peter Ozsv\'ath}
\address{Department of
Mathematics, Columbia University, New York 10025 \newline
\indent{Institute for Advanced Study, Princeton, New Jersey 08540} \newline
\indent{\emailfont{petero@math.columbia.edu}}}
\thanks{PSO was partially supported by NSF grant numbers DMS-0234311,
DMS-0111298,
and FRG-0244663}

\author[Zolt{\'a}n Szab{\'o}]{Zolt{\'a}n Szab{\'o}} 
\address{Department of
Mathematics, Princeton University, New Jersey 08544 \newline
\indent{\emailfont{szabo@math.princeton.edu}}}}
\thanks{ZSz was partially supported by NSF grant numbers DMS-0107792
and FRG-0244663, and a Packard Fellowship.}
 
\begin{document}

\maketitle
\section{Introduction}

Gromov's theory of pseudo-holomorphic disks~\cite{Gromov} has
wide-reaching consequences in symplectic geometry and low-dimensional
topology. Our aim here is to describe certain invariants for
low-dimensional manifolds built on this theory.

 The invariants we describe here associate a graded Abelian group to
each closed, oriented three-manifold $Y$, the Heegaard Floer homology
of $Y$. These invariants also have a four-dimensional counterpart,
which associates to each smooth cobordisms between two such
three-manifolds, a map between the corresponding Floer homology
groups. In another direction, there is a variant which gives rise to
an invariant of knots in $Y$.

\subsection{Some background on Floer homology}

To place Heegaard Floer homology into a wider context, we begin
with Casson's invariant. Starting with a Heegaard decomposition of an
integer homology three-sphere $Y$, Casson constructs a numerical invariant
which roughly speaking gives an obstruction to disjoining the $SU(2)$
character varieties of the two handlebodies inside the character
variety for the Heegaard surface $\Sigma$, c.f.~\cite{Casson}, \cite{Saveliev}.

During the time when Casson introduced his invariants to
three-dimensional topology, smooth four-dimensional
topology was being revolutionized by the work of
Donaldson~\cite{Donaldson}, who showed that the moduli spaces of
solutions to certain non-linear, elliptic PDEs -- gauge theory
equations which were first written down by physicists -- revealed a
great deal about the underlying smooth four-manifold topology.
Indeed, he constructed certain diffeomorphism invariants, called
Donaldson polynomials, defined by counting (in a suitable sense)
solutions to these PDEs, the anti-self-dual Yang-Mills
equations~\cite{Chambers}, \cite{DonaldsonPolynomials},
\cite{DonKron}, \cite{FriedmanMorgan}, \cite{KMPolyStruct}.

It was proved by Taubes in~\cite{TaubesCasson} that Casson's invariant
admits a gauge-theoretic interpretation. This interpretation was
carried further by Floer~\cite{InstantonFloer}, who constructed a
homology theory whose Euler characteristic is Casson's invariant.
The  construction of this instanton Floer homology
proceeds by defining a chain complex whose generators are equivalence
classes of flat $SU(2)$ connections over $Y$ (or, more precisely, a
suitably perturbed notion of flat connections, as required for
transversality), and whose differentials count solutions to the
anti-self-dual Yang-Mills equations. In fact, Floer's instanton
homology quickly became a central tool in the calculation of
Donaldson's invariants, see for example~\cite{MMR}, \cite{FSsfs},
\cite{DonaldsonBook}. More specifically, under suitable conditions,
the Donaldson invariant of a four-manifold $X$ separated along a
three-manifold $Y$ could be viewed as a pairing, in the Floer homology
of $Y$, of relative Donaldson invariants coming from the two sides.

Floer's construction seemed closely related to an earlier construction
Floer gave in the context of Hamiltonian dynamics, known as
``Lagrangian Floer homology''~\cite{FloerLag}. That theory -- which is
also very closely related to Gromov's invariants for symplectic
manifolds, c.f.~\cite{Gromov} -- associates to a symplectic manifold
$V$, equipped with a pair of Lagrangian submanifolds $L_0$ and $L_1$
(that in generic position, and satisfy certain topological
restrictions), a homology theory whose Euler characteristic is the
algebraic intersection number of $L_0$ and $L_1$, but which gives a
refined symplectic obstruction to disjoining the Lagrangians through
exact Hamiltonian isotopies. More specifically, the generators for
this chain complex are intersection points for $L_0$ and $L_1$, and
its differentials count holomorphic Whitney disks which interpolate
between these intersection points, see also~\cite{FOOO}.

The close parallel between Floer's two constructions, which take us
back to Casson's original picture, were further explored by
Atiyah~\cite{AtiyahFloer}. Atiyah conjectured that Floer's instanton
theory coincides with a suitably version of Floer's Lagrangian theory,
where one considers the $SU(2)$ character variety of $\Sigma$ as the
ambient symplectic manifold, equipped with the Lagrangian submanifolds
which are the character varieties of the two handlebodies. The
Atiyah-Floer conjecture remains open to this day. For related results,
see~\cite{DostoglouSalamon}, \cite{AFSalamon},
\cite{Wehrheim}.

In 1994, there was another drastic turn of events in gauge theory and
its interaction with smooth four-manifold topology, namely, the
introduction of a new set of equations coming from physics, the
Seiberg-Witten monopole equations~\cite{Witten}.  These are a novel
system of non-linear, elliptic, first-order equations which one can
associate to a smooth four-manifold equipped with a Riemannian metric.
Just as the Yang-Mills equations lead to Donaldson polynomials, the
Seiberg-Witten equations lead to another smooth four-manifold
invariant, the Seiberg-Witten invariant, c.f.~\cite{Witten},
\cite{Morgan},
\cite{DonaldsonSurvey}, \cite{KMthom}, \cite{TaubesBook}. 
These two theories seem very closely related. In fact, Witten
conjectured a precise relationship between the two four-manifold
invariants, see~\cite{Witten},
\cite{FeehanLeness}, see also~\cite{KMpropP}. Moreover, many of the formal aspects of Donaldson
theory have their analogues in Seiberg-Witten theory. In particular,
it was natural to expect a similar relationship between their
three-dimensional counterparts~\cite{MengTaubes}, \cite{MOY},
\cite{MarcolliWang}.

But now, a question arises in studying the three-dimensional theory:
what is the geometric picture playing the role of character varieties
in this new context? In attempting to formulate an answer to this
question, we came upon a construction which has, as its starting point
a Heegaard diagram $(\Sigma,\alphas,\betas)$ for a three-manifold
$Y$~\cite{HolDisk}.  That is, $\Sigma$ is an oriented surface of genus
$g$, and $\alphas=\{\alpha_1,...\alpha_g\}$ and
$\betas=\{\beta_1,...\beta_g\}$ are a pair of $g$-tuples of embedded,
homologically linearly independent, mutually disjoint, closed
curves. Thus, the $\alphas$ and $\betas$ specify a pair of
handlebodies $U_\alpha$ and $U_\beta$ which bound $\Sigma$, so that
$Y\cong U_\alpha\cup_{\Sigma} U_\beta$.  Note that any oriented,
closed three-manifold can be described by a Heegaard diagram. We
associate to $\Sigma$ its $g$-fold symmetric product $\Sym^g(\Sigma)$,
the space of unordered $g$-tuples of points in $\Sigma$. This space is
equipped with a pair of $g$-dimensional tori
\begin{eqnarray*}
\Ta=\alpha_1\times...\times \alpha_g&{\text{and}}&
\Tb=\beta_1\times...\times\beta_g.
\end{eqnarray*}
The most naive numerical invariant in this context -- the oriented
intersection number of $\Ta$ and $\Tb$ -- depends only on
$H_1(Y;\Z)$. However, by using the holomorphic disk techniques of
Lagrangian Floer homology, we obtain a non-trivial invariant for
three-manifolds, $\HFa(Y)$, whose Euler characteristic is this
intersection number. In fact, there are some additional elaborations
of this construction which give other variants of Heegaard Floer
homology (denoted $\HFm$, $\HFinf$, and $\HFp$, discussed below).

This geometric construction gives rise to invariants whose definition
is quite  different in flavor than its gauge-theoretic predecessors.  And
yet, it is natural to conjecture that certain variants give the same
information as Seiberg-Witten theory. This conjecture, in
turn, can be viewed as an analogue of the Atiyah-Floer conjecture in
the Seiberg-Witten context. With this said, it is also fruitful to
study Heegaard Floer homology and its structure independently from its
gauge-theoretical origins.

\subsection{Structure of this paper}

Our aim in this article is to give a leisurely introduction to
Heegaard Floer homology. We begin by recalling some of the details of
the construction in Section~\ref{sec:Construction}. In
Section~\ref{sec:BasProp}, we describe some of the properties. Broader
summaries can be found in some of our other papers
(c.f.~\cite{HolDiskTwo}, \cite{AbsGraded}). In Section~\ref{sec:Knots}
we describe in further detail the relationship between Heegaard Floer
homology and knots, c.f.~\cite{HolDiskKnots}, \cite{AltKnots} and also
the work of Rasmussen~\cite{Rasmussen}, \cite{RasmussenThesis}.  We
conclude in Section~\ref{sec:Questions} with some problems and
questions raised by these investigations.

We have not attempted to give a full account of the state of Heegaard
Floer homology. In particular, we have said very little about the
four-manifold invariants. We do not discuss here the Dehn surgery
characterization of the unknot which follows from properties of Floer
homology (see Corollary~\ref{GenusBounds:cor:KMOSz} of~\cite{GenusBounds},
and also~\cite{KMOS} for the original proof using Seiberg-Witten
monopole Floer homology; compare also~\cite{GordonLueckeI},
\cite{CGLS}, \cite{GabaiKnots}).  Another topic to which we have paid
only fleeting attention is the close relationship between Heegaard
Floer homology and contact geometry~\cite{HolDiskContact}. As a
result, we do not have the opportunity to describe the recent results
of Lisca and Stipsicz in contact geometry which result from this
interplay, see for example~\cite{LiscaStipsicz}.

\subsection{Further remarks}

The conjectured relationship between Heegaard Floer homology and
Seiberg-Witten theory can be put on a more precise footing with the
help of some more recent developments in gauge theory. For example,
Kronheimer and Mrowka~\cite{KMbook} have given a complete construction
of a Seiberg-Witten-Floer package, which associates to each closed,
oriented three-manifold a triple of Floer homology groups
$\HMto(Y)$, $\HMred(Y)$, and $\HMfrom(Y)$ which are functorial under
cobordisms between three-manifolds. They conjecture that the three
functors in this sequence are isomorphic to (suitable completions of)
$\HFm(Y)$, $\HFinf(Y)$ and $\HFp(Y)$ respectively. 
A different approach is taken in papers by Manolescu and
Kronheimer, c.f.~\cite{Manolescu}, \cite{ManolescuKronheimer}.

It should also be pointed out that a different approach to
understanding gauge theory from a geometrical point of view has been
adopted by Taubes~\cite{TaubesApproach}, building on his fundamental
earlier work relating the Seiberg-Witten and Gromov invariants of
symplectic four-manifolds~\cite{TaubesSympI}, \cite{TaubesSympII},
\cite{TaubesBook}.

\section{The construction}
\label{sec:Construction}

We recall the construction of Heegaard Floer homology. In
Subsection~\ref{subsec:RatSpheres}, we explain the construction for
rational homology three-spheres (i.e. those three-manifolds whose
first Betti number vanishes). In Subsection~\ref{subsec:Example},
we give an example which illustrates
some of the subtleties involved in the Floer complex.
In Subsection~\ref{subsec:BigMans}, we 
outline how the construction can be generalized for arbitrary closed,
oriented three-manifolds. In Subsection~\ref{subsec:Maps}, we sketch
the construction of the maps induced by cobordisms, and in
Subsection~\ref{subsec:Knots} we give preliminaries on the
construction of the invariants for knots in $S^3$.  The material in
Sections~\ref{subsec:RatSpheres}-\ref{subsec:BigMans} is derived
from~\cite{HolDisk}. The material from Subsection~\ref{subsec:Maps} is
an account of the material starting in
Section~\ref{HolDisk:sec:HolTriangles} if~\cite{HolDisk} and continued
in~\cite{HolDiskFour}.  The material from
Subsection~\ref{subsec:Knots} can be found in~\cite{HolDiskKnots}, see
also~\cite{RasmussenThesis}.

\subsection{Heegaard Floer homology for rational homology three-spheres}
\label{subsec:RatSpheres}

A {\em genus $g$ handlebody} $U$ is the three-manifold-with-boundary
obtained by attaching $g$ one-handles to a zero-handle. More
informally, a genus $g$ handlebody is homeomorphic to a regular
neighborhood of a bouquet of $g$ circles in $\R^3$.  The boundary of
$U$ is a two-manifold with genus $g$. If $Y$ is any oriented
three-manifold, for some $g$ we can write $Y$ as a union of two genus
$g$ handlebodies $U_{0}$ and $U_{1}$, glued together along their
boundary. A natural way of thinking about Heegaard decompositions is
to consider self-indexing Morse functions $$f\colon Y\longrightarrow
[0,3]$$ with one index $0$ critical point, one index three critical
point, and $g$ index one (hence also index two) critical points. The
space $U_{0}$, then, is the preimage of the interval $[0,3/2]$ (with
boundary the preimage of $3/2$); while $U_{1}$ is the preimage of
$[3/2,3]$. We orient $\Sigma$ as the boundary of $U_0$.

Heegaard decompositions give rise to a combinatorial description of
three-manifolds. Specifically, 
let $\Sigma$ be a closed, oriented surface of genus $g$. A {\em set of
attaching circles} for $\Sigma$ is a $g$-tuple of homologically
linearly independent, pairwise disjoint, embedded curves
$\gammas=\{\gamma_1,...,\gamma_g\}$. A {\em Heegaard diagram} is a
triple consisting of $(\Sigma,\alphas,\betas)$, where $\alphas$ and
$\betas$ are both complete sets of attaching circles for $\Sigma$.
From the Morse-theoretic point of view, the points in $\alphas$ can be
thought of as the points on $\Sigma$ which flow out of the index one
critical points (with respect to a suitable metric on $Y$), and the
points in $\betas$ are points in $\Sigma$ which flow into the index
two critical points.

In the opposite direction, a set of attaching circles for $\Sigma$
specifies a handlebody which bounds $\Sigma$, and hence a Heegaard
diagram specifies an oriented three-manifold $Y$.  It is a classical
theorem of Singer~\cite{Singer} that every closed, oriented
three-manifold $Y$ admits a Heegaard diagram, and if two Heegaard
diagrams describe the same three-manifold, then they can be connected
by a sequence of moves of the following type:
\begin{itemize}
    \item {\em isotopies:} replace $\alpha_{i}$ by a curve
    $\alpha_{i}'$ which is isotopic through isotopies which are
    disjoint from the other $\alpha_{j}'$ ($j\neq i$); or, the same
    moves amongst the $\betas$ \item {\em handleslides:} replace
    $\alpha_{i}$ by $\alpha_{i}'$, which is a curve with the property
    that $\alpha_{i}\cup \alpha_{i}'\cup \alpha_{j}$ bound a pair of
    pants which is disjoint from the remaining $\alpha_{k}$ ($k\neq
    i,j$); or, the same moves amongst the $\betas$ \item {\em
    stabilizations/destabilizations:} A stabilization replaces
    $\Sigma$ by its connected sum with a genus one surface
    $\Sigma'=\Sigma\# E$, and replaces 
    $\{\alpha_1,...,\alpha_g\}$ and $\{\beta_1,...,\beta_g\}$ by
    $\{\alpha_{1},\ldots,\alpha_{g+1}\}$ and
    $\{\beta_{1},\ldots,\beta_{g+1}\}$ respectively, where
    $\alpha_{g+1}$ and $\beta_{g+1}$ are a pair of curves supported in
    $E$, meeting transversally in a single point.
\end{itemize}

Our goal is to associate a group to each Heegaard diagram,
which is unchanged by the above three operations, and hence 
an invariant of the underlying three-manifold.

To this end, we will use a variant of Floer homology in the $g$-fold
symmetric product of a genus $g$ Heegaard surface $\Sigma$,
relative to the pair of 
totally real subspaces
\begin{eqnarray*}
\Ta=\alpha_1\times ...\times \alpha_g
&{\text{and}}&
\Tb=\beta_1\times ...\times \beta_g.
\end{eqnarray*}
That is to say, we define a chain complex generated by intersection
points between $\Ta$ and $\Tb$, and whose boundary maps count
pseudo-holomorphic disks in $\Sym^g(\Sigma)$.  Again, it is useful to
bear in mind the Morse-theoretic interpretation: an intersection point
$\x$ between $\Ta$ and $\Tb$ can be viewed as a $g$-tuple of gradient
flow-lines which connect all the index two and index one critical
points. We denote the corresponding one-chain in $Y$ by $\gamma_\x$
and call it a {\em simultaneous trajectory}.

Before studying the space of pseudo-holomorphic Whitney disks, 
we turn out attention to the algebraic topology of Whitney disks.
Specifically, we consider the unit disk
$\Strip$ in $\C$, and let $e_1\subset \partial \Strip$
denote the arc where $\Real(z)\geq 0$, and $e_2\subset \partial \Strip$
denote the arc where $\Real(z)\leq 0$.
Let $\pi_{2}(\x,\y)$
denote the set of homotopy classes of Whitney disks, i.e. maps
$$\left\{u\colon \Strip\longrightarrow
\Sym^g(\Sigma)\Bigg|
\begin{array}{l}
u(-i)=\x, u(i)=\y \\
u(e_1)\subset \Ta, u(e_2)\subset \Tb
\end{array}\right\}.$$                                                         

Fix intersection points $\x,\y\in\Ta\cap \Tb$. There is an obvious
obstruction to the existence of a Whitney disk which lives in
$H_1(Y;\Z)$.  It is obtained as follows. Given $\x$ and $\y$, consider
the corresponding simultaneous trajectories $\gamma_\x$ and
$\gamma_\y$.  The difference $\gamma_\x-\gamma_\y$ gives a closed loop
in $Y$, whose homology class is trivial if there exists a Whitney disk
connecting $\x$ and $\y$. (Indeed, the condition is sufficient when
$g>1$.)

What this shows is that the space of intersection points between 
$\Ta$ and $\Tb$ naturally fall into equivalence classes labeled by 
elements in an affine space over $H_{1}(Y;\Z)$. 

There is another very familiar such affine space over $H_{1}(Y;\Z)$:
the space of $\SpinC$ structures over $Y$. Following
Turaev~\cite{Turaev}, one can think of this space concretely as the
space of equivalence classes of vector fields. Specifically, we say
that two vector fields $v$ and $v'$ over $Y$ are {\em homologous} if they
agree outside a Euclidean ball in $Y$. The space of homology classes
of vector fields (which in turn are identified with the more standard
definitions of $\SpinC$ structures, see~\cite{Turaev},
\cite{GompfStipsicz},
\cite{KMcontact}) is also an affine space over $H_{1}(Y;\Z)$.

In order to link these two concepts, we fix a basepoint
$z\in\Sigma-\alpha_{1}-\ldots-\alpha_{g}-\beta_{1}-\ldots-\beta_{g}$. Thinking
of the Heegaard decomposition as a Morse function as described
earlier, the base point $z$ describes a flow from the index zero to
the index three critical point (for generic metric on $Y$). Now, each
tuple $\x\in\Ta\cap\Tb$ specifies a $g$-tuple of connecting flows
between the index one and index two critical points. Modifying the
gradient vector field $\grad f$ in a tubular neighborhood of these
$g+1$ flow-lines so that it does not vanish there, we obtain a nowhere
vanishing vector field over $Y$, whose homology class gives us a
$\SpinC$ structure, depending on $\x$ and $z$. This gives rise to an
assignment $$\spinc_{z}\colon \Ta\cap\Tb\longrightarrow
\SpinC(Y).$$
It follows from the previous discussion that $\x$ and $\y$ can be
connected by a Whitney disk if and only if
$\spinc_z(\x)=\spinc_z(\y)$.

Having answered the existence problem for
Whitney disks, we turn to questions of 
its uniqueness (up to homotopy). For this, we turn once again to 
our fixed base-point 
$z\in\Sigma-\alpha_{1}-\ldots-\alpha_{g}-\beta_{1}-\ldots-\beta_{g}$. 
Given a Whitney disk $u$ connecting $\x$ with $\y$, we can consider 
the intersection number between $u$ and the submanifold 
$\{z\}\times\Sym^{g-1}(\Sigma)\subset \Sym^{g}(\Sigma)$. This descends 
to homotopy classes to give a map 
$$n_{z}\colon \pi_{2}(\x,\y)\longrightarrow \Z$$
which is additive in the sense that
$n_{z}(\phi*\psi)=n_{z}(\phi)+n_{z}(\psi)$,
where here $*$ denotes the natural juxtaposition operation. 
The multiplicity $n_{z}$ can be modified by splicing in 
a copy of the two-sphere which generates 
$\pi_{2}(\Sym^{g}(\Sigma))\cong \Z$ (when $g>2$).
Moreover, when $Y$ is a rational 
homology sphere and $g>2$,  $n_{z}(\phi)$ uniquely 
determines the homotopy class of $\phi$  -- more generally, 
we have an identification $\pi_{2}(\x,\y)\cong \Z\oplus H^{1}(Y;\Z)$. 
For the present discussion, though, we focus attention to the rational 
homology sphere case.

The most naive application of Floer's theory would then give 
a $\Zmod{2}$-graded theory. However, 
the calculation of $\pi_{2}(\x,\y)$ suggests that if we count each 
intersection number infinitely many times, we obtain a relatively $\Z$-graded 
theory. Specifically, for a fixed $\SpinC$ structure $\spinc$ over 
$Y$,
let the set
${\mathfrak S}\subset \Ta\cap \Tb$
consist of intersection points which induce $\spinc$, 
with respect to the fixed base-point $z$. Now 
we can consider the Abelian group $\CFinf(Y,\spinc)$ freely generated 
by the set of pairs 
$[\x,i]\in{\mathfrak S}\times \Z$. We can give this space a natural 
relative
$\Z$-grading, by
\begin{equation}
    \label{eq:DefGrading}
    \gr([\x,i],[\y,j])=\Mas(\phi)-2(i-j+n_z(\phi)),
\end{equation}
where here $\phi$ is any homotopy class of Whitney disks which
connects $\x$ and $\y$, and $\Mas(\phi)$ denotes the Maslov index of
$\phi$, that is, the expected dimension of the moduli space of
pseudo-holomorphic representatives of $\phi$, see
also~\cite{RobbinSalamon}. If $S$ denotes the generator of
$\pi_2(\Sym^g(\Sigma))\cong \Z$ (when $g>2$), then
$\Mas(\phi*S)=\Mas(\phi)+2$ while 
$n_z(\phi*S)=n_z(\phi)+1$, and hence $\gr$ is
independent of the choice of $\phi$. (It is easy to see that $\gr$ extends
also to the cases where $g\leq 2$.)

Our aim will be to count pseudo-holomorphic disks.  For this to make
sense, we need to have a sufficiently generic situation, so that the
moduli spaces are cut out transversally and, in particular,
\begin{equation}
\label{eq:Transversality}
\dim\ModFlow(\phi)=\Mas(\phi).
\end{equation}
To achieve this, we need to introduce a suitable perturbation of the
notion of pseudo-holomorphic disks, see for
instance~Section~\ref{HolDisk:sec:Analysis} of~\cite{HolDisk}, see
also ~\cite{FloerHoferSalamon}, \cite{FOOO}. Specifically, for such a
perturbation, we can arrange Equation~\eqref{eq:Transversality} to
hold for all homotopy classes $\phi$ with $\Mas(\phi)\leq 2$.

Indeed, since there is a one-parameter family of
holomorphic automorphisms of the disks which preserve $\pm i$ and the
boundary arcs $e_{1}$ and $e_{2}$, the moduli space $\ModFlow(\phi)$
admits a free action by $\R$, provided that $\phi$ is
non-trivial.
In particular, if $\phi$ has $\Mas(\phi)=1$, then
$\ModFlow(\phi)/\R$ is a zero-dimensional manifold.

We then define a boundary map on $\CFinf(Y,\spinc)$ by the formula
\begin{equation}
    \label{eq:DefBoundary}
    \partial^\infty [\x,i] = 
	\sum_{\y}\sum_{\{\phi\in\pi_2(\x,\y)|\Mas(\phi)=1\}}
	\#\left(
    \frac{\ModFlow(\phi)}{\R}\right)[\y,i-n_z(\phi)].
\end{equation}
Here, $\#\left(\frac{\ModFlow(\phi)}{\R}\right)$ can be thought of as
either a  count modulo $2$ of the number of points in the moduli space,
in the case where we consider Floer homology with coefficients in
$\Zmod{2}$ or, in a more general case, to be an appropriately
signed count of the number of points in the moduli space. 
For a discussion on signs, see~\cite{HolDisk}, see also~\cite{FloerHofer},
\cite{FOOO}.

By analyzing Gromov limits of pseudo-holomorphic disks,
one can prove that $(\partial^\infty)^2=0$, i.e.
that $\CFinf(Y,\spinc)$ is a chain complex. 

It is easy to see that if a given homotopy class $\phi$ 
contains a holomorphic representative, then its intersection number 
$n_{z}(\phi)$ is non-negative. This observation ensures that the 
subset $\CFm(Y,\spinc)\subset \CFinf(Y,\spinc)$ generated by $[\x,i]$ 
with $i<0$ is a subcomplex. It is also interesting to consider the
quotient complex
$\CFp(Y,\spinc)$ (which we can think of as generated by pairs $[\x,i]$ 
with $i\geq 0$). Note that all three complexes can be thought of
as $\Z[U]$-modules, where 
$$U\cm [\x,i]=[\x,i-1];$$
i.e. multiplication by $U$ lowers grading by two.
Similarly, we can define a complex $\CFa(Y,\spinc)$ which 
is generated by the kernel of the $U$-action on $\CFp(Y,\spinc)$.
We can think of $\CFa(Y)$ directly as generated by intersection points 
between $\Ta$ and $\Tb$, endowed with the differential
\begin{equation}
\label{eq:Differentiala}
{\widehat \partial} \x = \sum_{\{\phi\in\pi_2(\x,\y)\big|\Mas(\phi)=1,
n_{z}(\phi)=0\}} \#\left(\frac{\ModFlow(\phi)}{\R}\right)\cm \y.
\end{equation}

We now define Floer homology theories $\HFm(Y,\spinc)$,
$\HFinf(Y,\spinc)$, $\HFp(Y,\spinc)$, $\HFa(Y,\spin)$,
which are the homologies of the chain complexes
$\CFm(Y,\spinc)$,
$\CFinf(Y,\spinc)$, $\CFp(Y,\spinc)$, and $\CFa(Y,\spin)$
respectively. Note that all of these groups are
$\Z[U]$ modules (where the $U$ action is
trivial on $\HFa(Y,\spinc)$).

The main result of~\cite{HolDisk}, then, is the following 
topological invariance of these theories:

\begin{theorem}
    \label{thm:TopInvar} The relatively $\Z$-graded theories
    $\HFm(Y,\spinc)$, $\HFinf(Y,\spinc)$, $\HFp(Y,\spinc)$, and
    $\HFa(Y,\spinc)$ are topological invariants of the underlying
    three-manifold $Y$ and its $\SpinC$ structure $\spinc$.
\end{theorem}

The content of the above result is that the invariants are independent 
of the various choices going into the definition of the homology 
theories. It can be broken up into parts, where one shows that the 
homology groups are identified as the Heegaard diagram undergoes 
the following changes:
\begin{enumerate}
    \item the complex structure over $\Sigma$ is varied
    \item the attaching circles are moved by isotopies (in the 
	complement of $z$)
    \item the attaching circles are moved by handle-slides (in the 
	complement of $z$)
    \item the Heegaard diagram is stabilized.
\end{enumerate}
The first step is a direct adaptation of the corresponding fact from
Lagrangian Floer theory (independence of the particular compatible
almost-complex structure). To see the second step, we observe that any
isotopy of the $\alphas$ and $\betas$ can be realized as a sequence of
exact Hamiltonian isotopies and metric changes over $\Sigma$. The
third step follows from naturality properties of the Floer homology
theories (using a holomorphic triangle construction which we return to
in Subsection~\ref{subsec:Maps}), and a direct calculation in a
special case (where handle-slides are made over a $g$-fold connected
sum of $S^{1}\times S^{2}$). The final step can be seen as an
invariance of the theory under a natural degeneration of the
$(g+1)$-fold symmetric product of the connected sum of $\Sigma$ with
$E$, as the connected sum neck is stretched, compare
also~\cite{LiRuan},
\cite{IonelParker}.

Although the study of holomorphic disks in general is a daunting task, 
holomorphic disks in symmetric products in a Riemann surface
admit a particularly nice interpretation in terms of the underlying Riemann 
surface: indeed, holomorphic disks in the $g$-fold symmetric product 
correspond to $g$-fold branched coverings of the disk by a Riemann 
surface-with-boundary, together with a holomorphic map from the 
Riemann-surface-with-boundary into $\Sigma$. This gives a geometric 
grasp of the objects under study, and hence, in many special cases, a 
way to calculate the boundary maps. In particular, it makes the model 
calculations used in the proof of handle-slide invariance mentioned 
above possible. We use it also freely in the following example.

\subsection{An example}
\label{subsec:Example}

We give a concrete example to illustrate some of the familiar
subtleties arising in the Floer complex. For simplicity, we stick to
the case of $\HFa(S^3)$.  Moreover, since we have made no attempt to
explain sign conventions, we consider the group with coefficients in
$\Zmod{2}$

\begin{figure}
\mbox{\vbox{\epsfbox{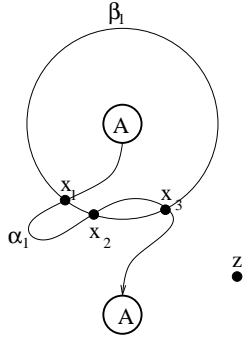}}}
\caption{\label{fig:SphereOne}
{\bf A genus one Heegaard diagram for $S^3$.}
In this diagram, the two circles labeled $A$ are to be identified, to obtain
a torus.}
\end{figure}

Of course, $S^3$ can be given a genus one Heegaard diagram, with two
attaching circles $\alpha_1$ and $\beta_1$, which meet in a unique
transverse intersection point. Correspondingly, the complex
$\CFa(S^3)$ in this case has a single generator, and there are no
differentials. Hence, $\HFa(S^3)\cong \Zmod{2}$. The diagram from
Figure~\ref{fig:SphereOne} is isotopic to this diagram, but now there
are three intersections between $\alpha_1$ and $\alpha_2$, $x_1$,
$x_2$, and $x_3$. By the Riemann mapping theorem, it is easy 
to see that ${\widehat \partial} x_1=x_2={\widehat \partial}x_3$.
Thus, $x_1+x_3$ generates $\HFa(S^3)$. Clearly, the chain complex
changed under the isotopy since the combinatorics of the new Heegaard diagram is
different (but, of course, its homology stayed the same).

\begin{figure}
\mbox{\vbox{\epsfbox{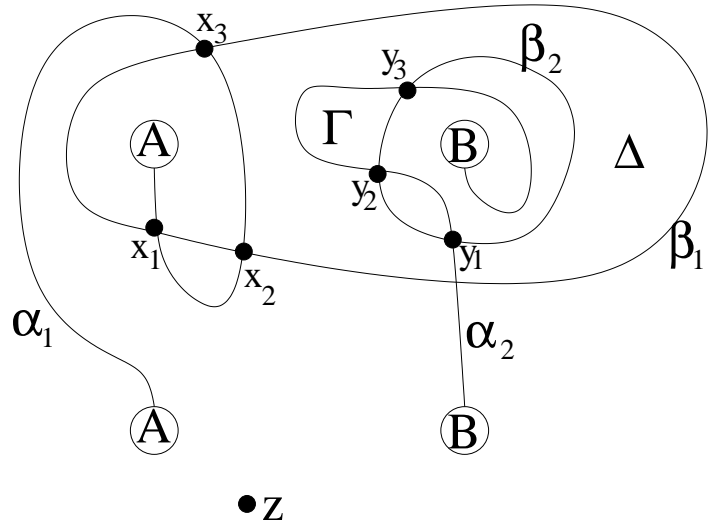}}}
\caption{\label{fig:GenusTwoDiagram}
{\bf A genus two Heegaard diagram for $S^3$.}
In this diagram, the two circles labeled $A$ are to be identified, as are
the circles labeled by $B$. The resulting surface $\Sigma$ of genus
two is divided into connected components by the union $\alpha_1\cup\alpha_2\cup\beta_1\cup\beta_2$. Let $\Delta$ be the component 
annular region indicated by taking the closure of the component indicated.}
\end{figure}

But the chain complex can change for reasons more subtle than
combinatorics.  Consider the Heegaard diagram for $S^3$ illustrated in
Figure~\ref{fig:GenusTwoDiagram}.

For this diagram, there are two different chain complexes, depending
on the choice of complex structure over $\Sigma$ (and the geometry of
the attaching circles). We sketch the argument below.

First, it is easy to see that there are nine generators, corresponding
to the points $x_i\times y_j\in\Sym^2(\Sigma)$ for $i,j=1,...,3$.
Again, by the Riemann mapping theorem applied to the region $\Gamma$,
there are holomorphic disks connecting $x_i\times y_3$ to $x_i\times
y_2$ for for all $i=1,...,3$. In a similar way, an inspection of
Figure~\ref{fig:GenusTwoDiagram} reveals disks connecting $x_1\times
y_j$ to $x_2\times y_j$ and $x_i\times y_1$ to $x_i\times
y_2$. However, the question of whether or not there is a holomorphic
disk in $\Sym^2(\Sigma)$ (with $n_z(\phi)=0$) connecting $x_3\times
y_i$ to $x_2\times y_j$ is dictated by the conformal structures of the
annuli in the diagram.

More precisely, consider the annular region $\Delta$ illustrated in
Figure~\ref{fig:GenusTwoDiagram}.  $\Delta$ has a uniformization as a
standard annulus with four points marked on its boundary,
corresponding to the points $x_1$, $x_2$, $y_2$, and $y_3$. Let $a$
denote the angle of the arc in the boundary connecting $x_1$ and $x_2$
which is the image of the corresponding segment in $\alpha_1$ under
this uniformization; let $b$ denote the angle of the arc in the
boundary connecting $y_2$ and $y_3$ which is the image of the
corresponding segment in $\alpha_2$ under the uniformization.  Now,
the question of whether there is a holomorphic disk in
$\Sym^2(\Sigma)$ connecting $x_3\times y_3$ to $x_2\times y_3$ admits
the following conformal reformulation.  Consider the one-parameter
family of conformal annuli with four marked boundary points obtained
from $\Delta\cup \Gamma$ by cutting a slit along $\alpha_2$ starting
at $y_3$. The four boundary points are the images of $x_3$, $x_2$, and
$y_3$ (counted twice) under a uniformization map. A four-times marked
annulus which admits an involution (interchanging the two
$\alpha$-arcs on the boundary) gives rise to a holomorphic disk
connecting $x_3\times y_3$ to $x_2\times y_3$. By analyzing the
conformal angles of the $\alpha$ arcs in this one-parameter family,
one can prove that the mod $2$ count of the holomorphic is $1$ iff
$a<b$.

\begin{figure}
\mbox{\vbox{\epsfbox{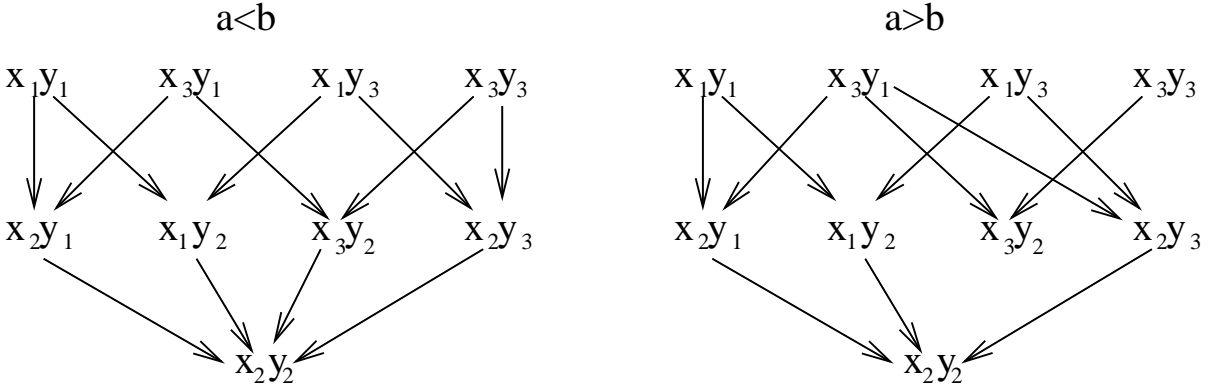}}}
\caption{\label{fig:Complexes}
{\bf Complexes for $\CFa(S^3)$ coming from the Heegaard diagram in
Figure~\ref{fig:GenusTwoDiagram}.}  The above two complexes can be
realized as $\CFa$ for a Heegaard diagram for $S^3$ illustrated in
Figure~\ref{fig:GenusTwoDiagram}, depending on the 
relations for the conformal parameter described in the text. Arrows
here indicate non-trivial differentials; e.g. for the complex on the left,
we have that ${\widehat \partial} x_1\times y_1 = x_2\times y_1 + x_1\times y_2$.}
\end{figure}

Proceeding in the like manner for the other homotopy classes, we see
that in the regime where $a<b$ resp. $a>b$, the chain complex $\CFa$
has differentials listed on the left resp. right in
Figure~\ref{fig:Complexes}. These two complexes are different, but of
course, they are chain homotopic.

\subsection{Algebra}
\label{subsec:Algebra}

The reason for this zoo of groups $\HFm$, $\HFinf$, $\HFp$, $\HFa$
can be traced to a simple algebraic reason:
$\CFm(Y,\spinc)$ (whose chain homotopy type is an invariant of $Y$) is
a finitely-generated chain complex of free $\Z[U]$-modules. All of the
other groups are obtained from this from canonical algebraic
operations. $\CFinf(Y,\spinc)$ is the ``localization''
$\CFm(Y,\spinc)\otimes_{\Z[U]} \Z[U,U^{-1}]$, $\CFp(Y,\spinc)$ is the
cokernel of the localization map, and $\CFa(Y,\spinc)$ is the quotient
$\CFm(Y,\spinc)/U \cm \CFm(Y,\spinc)$.

Correspondingly, the  various Floer homology groups are related by natural
long exact sequences
\begin{equation}
\label{eq:ExactSeqs}
\begin{CD}
...@>>>\HFm(Y,\spinc)@>{i}>>\HFinf(Y,\spinc) @>{\pi}>>\HFp(Y,\spinc)@>{\delta}>>...  \\
...@>>>\HFa(Y,\spinc)@>{j}>>\HFp(Y,\spinc) @>{U}>>\HFp(Y,\spinc)@>>>...,
\end{CD}
\end{equation}
and a more precise version of Theorem~\ref{thm:TopInvar} states that
both of the above diagrams are topological invariants of $Y$. The
interrelationships between these groups is essential in the study of
four-manifold invariants, as we shall see.

We can form another topological invariant, $\HFpRed(Y,\spinc)$, which
is the cokernel of $\pi$ appearing in Diagram~\eqref{eq:ExactSeqs}.

\subsection{Manifolds with $b_{1}(Y)>0$}
\label{subsec:BigMans}

When $b_{1}(Y)>0$, there are a number of additional technical issues
which arise
in the definition of Heegaard Floer homology. The crux of the matter
is that there are homotopically non-trivial cylinders
connecting $\Ta$ and $\Tb$. Specifically, given any point $\x$, we
have a subgroup of $\pi_{2}(\x,\x)$ consisting of classes of Whitney
disks $\phi$ with $n_{z}(\phi)=0$. This group, the group of ``periodic
classes,'' is naturally identified with the cohomology group
$H^{1}(Y;\Z)$, and hence (provided $g>2$) $\pi_{2}(\x,\y)\cong \Z\oplus
H^{2}(Y;\Z)$. In particular, there are infinitely many homotopy
classes of Whitney disks with a fixed multiplicity at a given point
$z$; thus, the coefficients appearing in
Equation~\eqref{eq:DefBoundary} might {\em a priori} be infinite.  One
way to remedy this situation is to work with special Heegaard diagrams
for $Y$. For example, in defining $\HFa$ and $\HFp$, 
one can use Heegaard diagrams with the property that
 for which each non-trivial periodic class has a negative
multiplicity at some
$z'\in\Sigma-\alpha_{1}-\ldots-\alpha_{g}-\beta_{1}-\ldots-\beta_{g}$.
Such diagrams are called {\em weakly admissible},
c.f. Section~\ref{HolDisk:subsec:Admissibility} of~\cite{HolDisk}
for a detailed account, and also for a discussion of the stronger
hypotheses needed for the construction of $\HFm$ and $\HFinf$.

Another related issue is that now, it is no longer true that the 
dimension of the space of holomorphic disks connecting $\x,\y$ 
depends only on the multiplicity at $z$. Specifically, given a 
one-dimensional cohomology class $\gamma\in H^{1}(Y;\Z)$, if 
$\phi\in\pi_{2}(\x,\y)$ and $\gamma*\phi$ denotes the new element of 
$\pi_{2}(\x,\y)$ obtained by letting the periodic class associated 
to $\gamma$ act on $\phi$, the Maslov classes of $\phi$ and 
$\gamma*\phi$ are related by the formula:
$$\Mas(\gamma*\phi)-\Mas(\phi)=\langle c_{1}(\spinc_{z}(\x))\cup \gamma,[Y]\rangle, $$
where the right-hand-side is, of course, calculated over the 
three-manifold $Y$. Letting $\delta(\spinc)$ be the greatest common divisor of 
the integers of the form $c_{1}(\spinc)\cup H^{1}(Y;\Z)$, 
the above discussion
shows that the grading defined in Equation~\eqref{eq:DefGrading} 
gives rise to a relatively $\Zmod{d(\spinc)}$-graded theory.

With this said, there is an analogue of Theorem~\ref{thm:TopInvar}:
when $b_{1}(Y)>0$, the homology theories $\HFp(Y,\spinc)$, 
$\HFm(Y;\spinc)$, and $\HFred(Y;\spinc)$ (as calculated for special 
Heegaard diagrams) are relatively $\Z/\delta(\spinc)$-graded
topological invariants. Note that, there are only finitely many $\SpinC$ 
structures $\spinc$ over $Y$ for which $\HFp(Y,\spinc)$ is non-zero.

\subsection{Maps induced by cobordisms}
\label{subsec:Maps}

Cobordisms between three-manifolds give rise to maps between their
Floer homology groups. The construction of these maps relies on the
holomorphic triangle construction from symplectic geometry,
c.f.~\cite{BraamDonaldson}, \cite{FOOO}.

A bridge between the symplectic geometry construction and the four-manifold 
picture can be given as follows. 

A {\em Heegaard triple-diagram of genus $g$} is an oriented two-manifold and three $g$-tuples
$\alphas$, $\betas$, and $\gammas$ which are complete sets of attaching
circles for handlebodies $U_{\alpha}$, $U_\beta$, and $U_{\gamma}$
respectively.
Let
$Y_{\alpha,\beta}=U_{\alpha}\cup U_{\beta}$,
$Y_{\beta,\gamma}=U_{\beta} \cup U_{\gamma}$, and
$Y_{\alpha,\gamma}=U_{\alpha} \cup U_{\gamma}$ denote the three
induced three-manifolds.
A Heegaard triple-diagram  naturally
specifies a cobordism $X_{\alpha,\beta,\gamma}$ between these
three-manifolds. The cobordism is constructed as follows.
                                                                                  
Let $\Delta$ denote the two-simplex, with vertices $v_{\alpha},
v_{\beta}, v_{\gamma}$ labeled clockwise, and let $e_{i}$ denote the
edge from $v_{j}$ to $v_{k}$, where $\{i,j,k\}=\{\alpha,\beta,\gamma\}$.
Then, we form the identification space
$$X_{\alpha,\beta,\gamma}=\frac{\left(\Delta\times \Sigma\right)
\coprod \left(e_\alpha \times U_\alpha\right) \coprod
\left(e_\beta\times U_\beta \right)\coprod
\left(e_\gamma\times U_\gamma \right)}
{ \left(e_{\alpha}\times \Sigma\right) \sim \left(e_\alpha \times
\partial U_\alpha\right), \left(e_{\beta}\times \Sigma\right) \sim
\left(e_\beta \times \partial U_\beta\right), \left(e_{\gamma}\times
\Sigma\right) \sim \left(e_\gamma \times \partial U_\gamma\right) }.
$$ Over the vertices of $\Delta$, this space has corners, which can be
naturally smoothed out to obtain a smooth, oriented, four-dimensional
cobordism between the three-manifolds $Y_{\alpha,\beta}$,
$Y_{\beta,\gamma}$, and $Y_{\alpha,\gamma}$ as claimed.

We will call the cobordism $X_{\alpha,\beta,\gamma}$ described above a
{\em pair of pants connecting $Y_{\alpha,\beta}$, $Y_{\beta,\gamma}$,
and $Y_{\alpha,\gamma}$}.  Note that if we give
$X_{\alpha,\beta,\gamma}$ its natural orientation, then $\partial
X_{\alpha,\beta,\gamma}=-Y_{\alpha,\beta}-Y_{\beta,\gamma}+Y_{\alpha,\gamma}$.
                                                                                
Fix $\x\in\Ta\cap \Tb$, $\y\in\Tb\cap \Tc$,
$\w\in\Ta\cap\Tc$. Consider the map $$u \colon \Delta \longrightarrow
\Sym^{g}(\Sigma) $$ with the boundary conditions that
$u(v_{\gamma})=\x$, $u(v_{\alpha})=\y$, and $u(v_{\beta})=\w$, and
$u(e_{\alpha})\subset \Ta$, $u(e_{\beta})\subset \Tb$,
$u(e_{\gamma})\subset
\Tc$. Such a map is called a {\em Whitney triangle connecting $\x$,
$\y$, and $\w$}. Two Whitney triangles are homotopic if the maps are
homotopic through maps which are all Whitney triangles. We let
$\pi_{2}(\x,\y,\w)$ denote the space of homotopy classes of Whitney
triangles connecting $\x$, $\y$, and $\w$.

Using a base-point
$z\in\Sigma-\alpha_1-...-\alpha_g-\beta_1-...-\beta_g-\gamma_1-...-\gamma_g$,
we obtain an intersection number
$$n_z\colon \pi_2(\x,\y,\w)\longrightarrow \Z.$$
If the space of homotopy classes of Whitney triangles
$\pi_2(\x,\y,\w)$ is non-empty, then it can be identified with
$\Z\oplus H_2(X_{\alpha,\beta,\gamma};\Z)$, in the case where $g>2$.

As explained in Section~\ref{HolDisk:sec:HolTriangles} of~\cite{HolDisk},
the choice of base-point $z\in\ThreeCurveComp$ gives rise to a map
$$\spinc_z\colon \pi_2(\x,\y,\w)\longrightarrow
\SpinC(X_{\alpha,\beta,\gamma}).$$ 
                                                                                                                                                          
A $\SpinC$ structure over $X$ gives rise to a map
$$
f^{\infty}(~\cdot~;\spinc)\colon
\CFinfty(Y_{\alpha,\beta},\spinc_{\alpha,\beta})
\otimes \CFinfty(Y_{\beta,\gamma},\spinc_{\beta,\gamma}) \longrightarrow
\CFinfty(Y_{\alpha,\gamma},\spinc_{\alpha,\gamma}) $$
by the formula:
\begin{equation}
\label{eq:DefTriangles}
f^{\infty}_{\alpha,\beta,\gamma}([\x,i]\otimes[\y,j];\spinc)
=\sum_{\w\in\Ta\cap\Tc}
\sum_{\{\psi\in\pi_2(\x,\y,\w)\big|\spinc_z(\psi)=\spinc,\Mas(\psi)=0\}}
\Big(\#\ModFlow(\psi)\Big)\cm {[\w,i+j-n_z(\psi)]}.
\end{equation}
Under suitable admissibility hypotheses on the Heegaard diagrams,
these sums are finite, c.f. Section~\ref{HolDisk:sec:HolTriangles}
of~\cite{HolDisk}. Indeed, there are induced maps on some of the other
variants of Floer homology, and again, we refer the interested reader
to that discussion for a more
detailed account.

Let $X$ be a smooth, connected, oriented four-manifold with boundary
given by $\partial X=-Y_0\cup Y_1$ where $Y_0$ and $Y_1$ are
connected, oriented three-manifolds. We call such a four-manifold a
cobordism from $Y_0$ to $Y_1$.  If $X$ is a cobordism from $Y_0$ to
$Y_1$, and $\spinc\in\SpinC(X)$ is a $\SpinC$ structure, then there is
a naturally induced map $$F^\infty_{X,\spinc}\colon
\HFinf(Y_0,\spinc_i)\longrightarrow \HFinf(Y_1,\spinc_i)$$ where here
$\spinc_i$ denotes the restriction of $\spinc$ to $Y_i$.  This map is
constructed as follows. First assume that $X$ is given as a collection
of two-handles. Then we claim that in the complement of the regular
neighborhood of a one-complex, $X$ can be realized as a pair-of-pants
cobordism, one of whose boundary components is $-Y_0$, the other which
is $Y_1$, and the third of which is a connected sum of copies of
$S^2\times S^1$. Next, pairing Floer homology classes coming from
$Y_0$ with a certain canonically associated Floer homology class on
the connected sum of $S^2\times S^1$, we obtain a map using the
holomorphic triangle construction as defined in
Equation~\eqref{eq:DefTriangles} to obtain a the desired map to
$\HFinf(Y_1)$. For the cases of one- and three-handles, the associated
maps are defined in a more formal manner. The fact that these maps are
independent (modulo an overall multiplication by $\pm 1$) of the many
choices which go into their construction is established
in~\cite{HolDiskFour}.

Indeed, variants of this construction can be extended to the following
situation (again, see~\cite{HolDiskFour}): if $X$ is a smooth,
oriented cobordism from $Y_0$ to $Y_1$, then there are induced maps
(of $\Z[U]$ modules) between the corresponding Heegaard Floer homology
groups, which make the squares in the following diagrams commutate:
\begin{equation}
\label{eq:NatTranses}
\begin{CD}
...@>>>\HFm(Y_0,\spinc_0)@>>>\HFinf(Y_0,\spinc_0) @>{\pi}>>\HFp(Y_0,\spinc_0)@>{\delta}>>...  \\
&& @V{F^-_{X,\spinc}}VV @V{F^\infty_{X,\spinc}}VV @V{F^+_{X,\spinc}}VV \\
...@>>>\HFm(Y_1,\spinc_1)@>>>\HFinf(Y_1,\spinc_1) @>{\pi}>>\HFp(Y_1,\spinc_1)@>{\delta}>>...  \\ \\
...@>>>\HFa(Y_0,\spinc_0)@>>>\HFp(Y_0,\spinc_0) @>{U}>>\HFp(Y_0,\spinc_0)@>>>..., \\
&& @V{{\widehat F}_{X,\spinc}}VV @V{F^+_{X,\spinc}}VV @V{F^+_{X,\spinc}}VV \\
...@>>>\HFa(Y_1,\spinc_1)@>>>\HFp(Y_1,\spinc_1) @>{U}>>\HFp(Y_1,\spinc_1)@>>>...,
\end{CD}
\end{equation}

Naturality of the maps induced by cobordisms can be phrased as follows.
Suppose that $W_0$ is a smooth cobordism from $Y_0$ to $Y_1$ and $W_1$
is a cobordism from $Y_1$ to $Y_2$, then for fixed $\SpinC$ structures
$\spinc_i$ over $W_i$ which agree over $Y_1$, we have that
$$\sum_{\{\spinc\in\SpinC(W_0\cup_{Y_1} W_1)\big| \spinc|_{W_i}=\spinc_i\}}
\Fc_{W_0\cup_{Y_1} W_1,\spinc}=
\Fc_{W_1,\spinc_1}\circ \Fc_{W_0,\spinc_0},$$
where here $\Fc=F^-$, $F^{\infty}$, $F^+$, ${\widehat F}$
(c.f. Theorem~\ref{HolDiskFour:thm:Composition} of~\cite{HolDiskFour}).

Sometimes, it is convenient to obtain topological invariants by
summing over all $\SpinC$ structures. To this end, we write, for example,
$$\HFp(Y)\cong \bigoplus_{\spinc\in\SpinC(Y)}\HFp(Y,\spinc).$$

It is convenient to have  a corresponding notion for cobordisms, only in that
case a  little more care must be taken.
For fixed $X$ and $\xi\in\HFp(Y_0,\spinc_0)$, we have that
$F^+_{X,\spinc}(\xi)=0$ for all but finitely many
$\spinc\in\SpinC(X)$, c.f. Theorem~\ref{HolDiskFour:thm:Finiteness}
of~\cite{HolDiskFour}, and hence there is a well-defined map
$$F^+_{X}\colon \HFp(Y_0)
\longrightarrow \HFp(Y_1),$$
defined by
$$F^+_{X}=\sum_{\spinc\in\SpinC(X)}F^+_{X,\spinc}.$$ 
(note that the same construction works for $\HFa$, but it
does  not work for $\HFm$, $\HFinf$: for a given
$\xi\in\HFinf(Y_0)$, there might be infinitely many different
$\spinc\in\SpinC(X)$ for which $F^\infty_{X,\spinc}(\xi)$ is
non-zero).

\subsection{Doubly-pointed Heegaard diagrams and knot invariants}
\label{subsec:Knots}

Additional basepoints give rise to additional filtrations on Floer homology.
These additional filtrations can be given topological interpretations. We 
consider the case of two basepoints.

Specifically, a Heegaard diagram $(\Sigma,\alphas,\betas)$ for $Y$
equipped with two basepoints $w$ and $z$ gives rise to a knot in $Y$
as follows.  We connect $w$ and $z$ by a curve $a$ in
$\Sigma-\alpha_1-...-\alpha_g$ and also by another curve $b$ in
$\Sigma-\beta_1-...-\beta_g$. By pushing $a$ and $b$ into $U_{\alpha}$
and $U_{\beta}$ respectively, we obtain a knot $K\subset Y$.  We call
the data $(\Sigma,\alphas,\betas,w,z)$ a {\em doubly-pointed Heegaard
diagram compatible with the knot $K\subset Y$}. Given a knot $K$ in
$Y$, one can always find such a Heegaard diagram.

This can be thought of from the following Morse-theoretic point of
view. Let $Y$ be an oriented three-manifold, equipped with a
Riemannian metric and a self-indexing Morse function $$f\colon
Y\longrightarrow [0,3]$$ with one index $0$ critical point, one index
three critical point, and $g$ index one (hence also index two)
critical points. The knot $K$ now is obtained from the union of the
two flows connecting the index $0$ to the index $3$ critical points
which pass through $w$ and $z$.  We call a Morse function as in the
above construction one which is compatible with $K$.  Note also that
an ordering of $w$ and $z$ is equivalent to an orientation on
$K$. However, the invariants we construct can be shown to be
independent of the orientation of $K$, see~\cite{HolDiskKnots},
\cite{RasmussenThesis}.

The simplest construction now is to consider a differential on $\Ta\cap\Tb$
defined analogously to Equation~\eqref{eq:Differentiala},
only now we count holomorphic disks for which $n_z(\phi)=n_w(\phi)=0$.
More generally, we use the reference point $w$
to construct  the Heegaard Floer complex for $Y$, and then use the additional
basepoint $z$ to induce a filtration on this complex.
We describe this in detail for the case of knots in $S^3$, and using  the
chain complex $\CFa(S^3)$.

There is a unique function  $\Filt\colon \Ta\cap\Tb
\longrightarrow \Z$ satisfying the relation
\begin{equation}
\label{eq:FilDiff}
\Filt(\x)-\Filt(\y)=n_{z}(\phi)-n_w(\phi),
\end{equation}
for any $\phi\in\pi_2(\x,\y)$, and the additional symmetry
$$\#\{\x\in\Ta\cap\Tb\big| \Filt(\x)=i\} \equiv
\#\{\x\in\Ta\cap\Tb\big| \Filt(\x)=-i\}\pmod{2}
$$ for all
$i$ (compare, more generally, Equation~\eqref{eq:EulerChar}). (Alternatively, a more intrinsic characterization can be given in
terms of relative $\SpinC$ structures on the knot complement.) 
Clearly, if $\y$ appears in ${\widehat\partial}(\x)$ 
with non-zero multiplicity, then the homotopy class
$\phi\in\pi_2(\x,\y)$ with $n_w(\phi)=0$ admits a holomorphic
representative, and hence $\Filt(\x)-\Filt(\y)\geq 0$. Thus, any
filtration satisfying Equation~\eqref{eq:FilDiff} induces a filtration
on the complex $\CFa(S^3)$, by the rule that
$\Filt(K,i)\subset \CFa(S^3)$ is the subcomplex generated by
$\x\in\Ta\cap\Tb$ with $\Filt(\x)\leq i$.

It is shown in~\cite{HolDiskKnots} and~\cite{RasmussenThesis} that the
chain homotopy type of this filtration is a knot invariant.  More
precisely, recall that a {\em filtered chain complex} is 
a chain complex $C$, together with a sequence of subcomplexes $X(C,i)$
indexed by $i\in\Z$, where $X(C,i)\subseteq X(C,i+1)\subset C$. 
Our filtered complexes are always bounded, meaning that for
all sufficiently large $i$, $X(C,-i)=0$ and $X(C,i)=C$. A {\em
filtered map} between chain complexes $\Phi\colon C \longrightarrow
C'$ is one whose restriction to $X(C,i)\subset C$ is contained in
$X(C',i)$.  Two filtered chain complexes are said to be have the same
{\em filtered chain type} if there are maps $f\colon C \longrightarrow
C'$ and $f'\colon C' \longrightarrow C$ and $H\colon C \longrightarrow
C'$ and $H'\colon C'\longrightarrow C'$, all four of which are
filtered maps, $f$ and $f'$ are chain maps, and also we have that
\begin{eqnarray*}
f\circ f' - \Id = \partial' \circ H' + H'\circ \partial' 
&{\text{and}}&
f'\circ f - \Id = \partial \circ H + H\circ \partial.
\end{eqnarray*}

The construction we mentioned earlier -- counting holomorphic disks
with $n_w(\phi)=n_z(\phi)=0$ can be thought of as the chain complex of
the associated graded object $$\bigoplus_{i}
\Filt(K,i)/\Filt(K,i-1).$$ The homology of this is also a knot
invariant. We return to properties of this invariant in
Section~\ref{sec:Knots}.

\section{Basic properties}
\label{sec:BasProp}

We outline here some of the basic properties of Heegaard Floer
homology, to give a flavor for its structure. We have not attempted to
summarize all of its properties; for additional properties,
see~\cite{HolDiskTwo}, \cite{AbsGraded}, \cite{HolDiskFour}.

We focus on material which is useful for calculations: an exact
sequence and rational gradings. We then turn briefly to properties of
the maps induced on $\HFinfty$, which have some important consequences
explained later, but they also shed light on the special role played
by $b_2^+(X)$ in Heegaard Floer homology. In
Section~\ref{subsec:Examples} we give a few sample calculations. In
Section~\ref{sec:IntFormBounds}, we describe one of the first
applications of the rational gradings: a constraint on the
intersection forms of four-manifolds which bound a given
three-manifold, compare the gauge-theoretic analogue of
Fr{\o}yshov~\cite{Froyshov}. Finally, in
Subsection~\ref{subsec:FourManInv}, we sketch how the maps induced by
cobordisms give rise to an interesting invariant of closed, smooth
four-manifolds $X$ with $b_2^+(X)>1$, which are conjectured to agree
with the Seiberg-Witten invariants, c.f.~\cite{Witten}.

\subsection{Long exact sequences}
\label{sec:LongExactSequence}

An important calculational device is provided by the surgery long
exact sequence. Long exact sequences of this type were first explored
by Floer in the context of instanton Floer homology~\cite{FloerTriangles},
\cite{BraamDonaldson}, see also~\cite{SeidelExactSeq}, \cite{KMOS}.

Heegaard Floer homology satisfies a surgery long exact sequence, which
we state presently.  Suppose that $M$ is a three-manifold with torus
boundary, and fix three simple, closed curves $\gamma_0$, $\gamma_1$,
and $\gamma_2$ in $\partial M$ with
\begin{equation}
\label{eq:TriadRelation}
\#(\gamma_0\cap \gamma_1)=\#(\gamma_1 \cap \gamma_2)
=\#(\gamma_2\cap \gamma_0) = -1
\end{equation} (where here the algebraic intersection
number is calculated in $\partial M$, oriented as the boundary of
$M$), so that $Y_0$ resp. $Y_1$ resp. $Y_2$ are obtained from $M$ by
attaching a solid torus along the boundary with meridian $\gamma_0$
resp. $\gamma_1$ resp. $\gamma_2$.

\begin{theorem}
\label{thm:ExactSeq}
Let $Y_0$, $Y_1$, and $Y_2$ be related as above. Then, 
there is a long exact sequence relating the Heegaard Floer homology groups:
$$
\begin{CD}
...@>>>\HFp(Y_0) @>>>\HFp(Y_1) @>>>\HFp(Y_2)@>>>...
\end{CD}
$$
\end{theorem}

The above theorem is proved in
Theorem~\ref{HolDiskTwo:thm:GeneralSurgery} of~\cite{HolDiskTwo}.
A variant for Seiberg-Witten monopole Floer homology, with coefficients 
in $\Zmod{2}$ is proved in~\cite{KMOS}. 

The maps in the long exact sequence have a four-dimensional
interpretation.  To this end, note that there are two-handle
cobordisms $W_i$ connecting $Y_i$ to $Y_{i+1}$ (where we view
$i\in\Zmod{3}$). When we work with Heegaard Floer homology over the
field $\Zmod{2}$, the map from $\HFp(Y_i)$ to $\HFp(Y_{i+1})$ in the
above exact sequence is the map induced by the corresponding cobordism
$W_{i}$, $F^+_{W_i}$ (i.e. obtained by summing the maps induced by all
$\SpinC$ structures over $W_i$). When working over $\Z$, though, one
must make additional choices of signs to ensure that exactness holds.

\subsection{Gradings}
\label{subsec:Gradings}

It is proved in Section~\ref{HolDiskTwo:sec:HFinfty} of~\cite{HolDiskTwo}
that if $Y$ is a rational homology three-sphere and $\spinc$
is any $\SpinC$ structure over it, then
$\HFinf(Y,\spinc)\cong \Z[U,U^{-1}]$, thus, this invariant is
not a very subtle invariant of three-manifolds. However, extra information
can still be gleaned from the interplay between $\HFinf$ and $\HFp$, 
with the help of some additional structure on Floer homology.

It is shown in~\cite{HolDiskFour} that when $Y$ is an oriented
rational homology three-sphere and $\spinc$ is a $\SpinC$ structure
over $Y$, the relative $\Z$ grading on the Heegaard Floer homology
described earlier can be lifted to an absolute $\Q$-grading.  This
gives $\HFc(Y,\spinc)$ is a $\Q$-graded module over the polynomial
algebra $\Z[U]$ (where here $\HFc(Y,\spinc)$ is any of
$\HFm(Y,\spinc)$, $\HFinf(Y,\spinc)$, $\HFp(Y,\spinc)$, or
$\HFa(Y,\spinc)$), $$\HFc(Y,\spinc)=\bigoplus_{d\in\Q}
\HFc_d(Y,\spinc),$$ where multiplication by $U$ lowers degree by
two. In each grading, $i\in\Q$, $\HFc_i(Y,\spinc)$ is a finitely
generated Abelian group.

The maps $i$, $\pi$, and $j$ in 
Diagram~\eqref{eq:ExactSeqs}
preserve this $\Q$-grading, and moreover, maps induced
by cobordisms $F^\circ_{X,\spinc}$ (again,
$F^\circ_{X,\spinc}$ denotes any of 
$F^-_{X,\spinc}$, $F^\infty_{X,\spinc}$,
$F^+_{X,\spinc}$ or ${\widehat F}_{X,\spinc}$)
respect the $\Q$-grading in the following sense. If
$Y_0$ and $Y_1$ are rational homology three-spheres, and 
$X$ is a cobordism from $Y_0$ to $Y_1$, with $\SpinC$ structure $\spinc$,
the map induced by the cobordism maps
$$F^{\circ}_{X,\spinc}\colon \HFc_d(Y_0,\spinc_0) \longrightarrow
\HFc_{d+\Delta}(Y_1,\spinc_1),$$
for
\begin{equation}
\label{eq:GradingShift}
\Delta = \frac{c_1(\spinc)^2-2\chi(X)-3\sigma(X)}{4},
\end{equation}
where here
$\chi(X)$ denotes the Euler characteristic of $X$, and
$\sigma(X)$ denotes its signature.  In fact
(c.f. Theorem~\ref{HolDiskFour:thm:AbsGrade} of~\cite{HolDiskFour})
the $\Q$ grading is uniquely characterized by the above property,
together with the fact that $d(S^3)=0$.

The image of $\pi$ determines a function $$d\colon \SpinC(Y)
\longrightarrow \Q$$ (the ``correction terms'' of~\cite{AbsGraded}) which associates to each $\SpinC$ structure the
minimal $\Q$-grading of any (non-zero) homogeneous element in
$\HFp(Y,\spinc)\otimes_\Z \Q$ in the image of $\pi$.  

Certain properties of the correction terms can be neatly summarized,
with the help of the following definitions.

The three-dimensional $\SpinC$ homology bordism group $\SpinCCobord$ is
the set of equivalence classes of 
pairs $(Y,\spinct)$ where $Y$ is a rational homology three-sphere,
and $\spinct$ is a $\SpinC$ structure over $Y$, and the equivalence
relation identifies $(Y_1,\spinct_1)\sim (Y_2,\spinct_2)$ if there is a
(connected, oriented, smooth) cobordism $W$ from $Y_1$ to $Y_2$ with
$H_i(W,\Q)=0$ for $i=1$ and $2$, which can be endowed with a $\SpinC$
structure $\spinc$ whose restrictions to $Y_1$ and $Y_2$ are $\spinct_1$
and $\spinct_2$ respectively. The connected sum operation endows this
set with the structure of an Abelian group (whose unit is $S^3$
endowed with its unique $\SpinC$ structure). 

There is a classical homomorphism $$\rho\colon
\theta^c\longrightarrow \Q/2\Z$$ (see for instance~\cite{APSII}),  defined
as follows. Consider a rational homology three-sphere $(Y,\spinct)$,
and let $X$ be any four-manifold equipped with a $\SpinC$ structure
$\spinc$ with $\partial X \cong Y$ and $\spinc|\partial X \cong
\spinct$. Then $$\rho(Y,\spinct)\equiv \frac{c_1(\spinc)^2-\sigma(X)}{4}
\pmod{2\Z}$$ where $\sigma(X)$ denotes the signature of the intersection
form of $X$.

It is shown in~\cite{AbsGraded} that 
the numerical invariant $d(Y,\spinct)$ descends to give a group homomorphism
$$d\colon \SpinCCobord\longrightarrow \Q$$
which is a lift of $\rho$. Moreover, 
$d$ is  invariant under conjugation; i.e.
$d(Y,\spinct)=d(Y,{\overline \spinct})$.

The rational gradings can be introduced for three-manifolds with $b_1(Y)>0$,
as well, only there one must restrict to $\SpinC$ structures
whose first Chern class is trivial. In this case, the gradings 
are fixed so that Equation~\eqref{eq:GradingShift} still holds. With these
conventions, for example, for a three-manifold with $H_1(Y;\Z)\cong \Z$,
the Heegaard Floer homologies of $\HFc(Y,\spinc)$ for the $\SpinC$ structure
with $c_1(\spinc)=0$ have a grading which takes its values in $\OneHalf+\Z$.

\subsection{Maps on $\HFinfty$}
\label{subsec:HFinfty}

As we have seen, for rational homology three-spheres, the structure
of $\HFinf$ is rather simple. There are corresponding statements
for the maps on $\HFinfty$ induced by cobordisms.

Indeed, if $W$ is a cobordism from $Y_1$ to $Y_2$ with $b_2^+(W)>0$,
the induced map $\Finf{W,\spinc}=0$ for any $\spinc\in\SpinC(W)$
(c.f. Lemma~\ref{HolDiskFour:lemma:BTwoPlusLemma}
of~\cite{HolDiskFour}). Moreover, if $W$ is a cobordism
from $Y_1$ to $Y_2$ (both of which are rational homology three-spheres),
and $W$ satisfies $b_2^+(W)=b_1(W)=0$, then $F^{\infty}_{W,\spinc}$
is an isomorphism, as proved in Propositions~\ref{AbsGraded:prop:ZeroSurgery} 
and~\ref{AbsGraded:prop:NegSurgery} of~\cite{AbsGraded}.

\subsection{Examples}
\label{subsec:Examples}

We begin with some algebraic notions for describing Heegaard Floer
homology groups. Let $\InjMod_{(d)}$ denote the graded $\Z[U]$-module
$\Z[U,U^{-1}]/\Z[U]$, graded so that the element $1$ has grading $d$.

A rational homology three-sphere $Y$ is called an {\em $L$-space} if
$\HFp(Y)$ has no torsion and the map from $\HFinfty(Y)$ to $\HFp(Y)$
is surjective. The Floer homology of an $L$-space can be uniquely
specified by its correction terms.  That is, if $Y$ is an $L$-space,
then $$\HFa(Y,\spinc)\cong
\Z_{(d(Y,\spinc))},$$ 
where here (and indeed throughout this subsection) 
the subscript denotes absolute grading,
and $$\HFp(Y,\spinc)\cong \InjMod_{(d)}.$$

By a direct inspection of the corresponding genus one Heegaard
diagrams, one can see that $S^3$ is an $L$-space. Indeed, by a similar
picture, all lens spaces are $L$-spaces.

The absolute $\Q$ grading can also be calculated for lens
spaces~\cite{AbsGraded}. For example, for $L(2,1)\cong \RP3$,
there are two $\SpinC$ structures $\spinc$ and $\spinc'$
with correction terms $1/4$ and $-1/4$ respectively.

The Brieskorn homology sphere $\Sigma(2,3,5)$ is also an $L$-space,
and it has $d(\Sigma(2,3,5))=2$. 

However, $\Sigma(2,3,7)$ is not an $L$-space. Its Heegaard Floer
homology is determined by 
\[\HFp(\Sigma(2,3,7))\cong 
\InjMod_{(0)} \oplus \Z_{(-1)}.\]

A combinatorial description of the Heegaard Floer homology of Brieskorn
spheres and some other plumbings can be found in ~\cite{SomePlumbs};
see also~\cite{AbsGraded}, \cite{Nemethi}, \cite{Rustamov}. 

\subsection{Intersection form bounds}
\label{sec:IntFormBounds}

The correction terms of a rational homology three-sphere $Y$ constrain
the intersection forms of smooth four-manifolds which bound $Y$,
according to the following result, which is analogous to a
gauge-theoretic result of Fr{\o}yshov~\cite{Froyshov}:

\begin{theorem}
  \label{thm:IntForms} Let $Y$ be a rational homology and $W$ be a
  smooth four-manifold which bounds $Y$ with negative-definite
  intersection form. Then, for each $\SpinC$ structure $\spinc$
  over $W$, we have that 
\begin{equation}
        \label{eq:Pos}
        c_1(\spinc)^2+b_2(W)\leq
        4 d(Y,\spinc|_{Y}).
\end{equation}
\end{theorem}

The above theorem gives strong restrictions on the intersection forms
of four-manifolds which bound a given three-manifold $Y$. In
particular, if $Y$ is an integral homology three-sphere, following a
standard argument from Seiberg-Witten theory, compare~\cite{Froyshov},
one can combine the above theorem with a number-theoretic result of
Elkies~\cite{Elkies} to show that if $Y$ can be realized as the
boundary of a smooth, negative-definite four-manifold, then $d(Y)\geq
0$; moreover if $d(Y)=0$, then if $X$ has negative-definite
intersection form, then it must be diagonalizable.

\subsection{Four-manifold invariants}
\label{subsec:FourManInv}

The invariants associated to cobordisms can be used to construct
an invariant for smooth, closed four-manifolds which is very similar in 
spirit to the Seiberg-Witten invariant for four-manifolds.
Indeed, all known calculations support the conjecture that the
two smooth four-manifold invariants agree.

Suppose that $X$ is a four-manifold with $b_2^+(X)>1$. We delete
four-ball neighborhoods of two points in $X$, and view the result as a
cobordism from $S^3$ to $S^3$, which we can further subdivide along a
separating hypersurface $N$ into a union $W_1\cup_{N} W_2$, with the following properties:
\begin{itemize}
\item $W_1$ is a cobordism from $S^3$ to $N$ with $b_2^+(W_1)>0$,
\item $W_2$ is a cobordism from $N$ to $S^3$ with $b_2^+(W_2)>0$,
\item restriction map $H^2(W_1\cup_N W_2) \longrightarrow 
 H^2(W_1)\oplus H^2(W_2)$ is injective.
\end{itemize}
Such a separating hypersurface is called an {\em admissible cut} for $X$.

Let $\HFmRed(Y)$ denote the kernel of the map $\HFm(Y)\longrightarrow
\HFinf(Y)$.  Of course, this is isomorphic to the group $\HFpRed(Y)$,
via an identification coming from the homomorphism $\delta$ from
Equation~\eqref{eq:NatTranses}.  Since $b_2^+(W_i)>0$, the maps on
$\HFinf$ induced by cobordisms are trivial (see
Lemma~\ref{HolDiskFour:lemma:BTwoPlusLemma} of~\cite{HolDiskFour}),
and in particular the image of the map $$\Fm{W_1,\spinc|W_1}\colon
\HFm(S^3)\longrightarrow
\HFm(N,\spinc|N)$$
lies in the kernel  
$\HFmred(N,\spinc|N)$ of the map $i$ (c.f. Diagram~\eqref{eq:ExactSeqs}. Moreover, the
map 
$$\Fp{W_2,\spinc|{W_2}}\colon \HFp(N,\spinc|{N})
\longrightarrow \HFp(S^3).$$ factors through
the projection of $\HFp(N,\spinc|N)$ to $\HFpred(N,\spinc|N)$ (the
cokernel of the map $\pi$ from Diagram~\eqref{eq:ExactSeqs}).  Thus,
we can define $$\Phi_{X,\spinc}\colon \HFm(S^3)
\longrightarrow 
\HFp(S^3)$$
to be the composite:
$$\Fp{W_2,\spinc|W_2}\circ \delta^{-1}\circ \Fm{W_1,\spinc|W_1},$$
where
$$\delta'\colon \HFpred(N,\spinc|N)\longrightarrow
\HFmred(N,\spinc|N)$$
is the natural isomorphism induced from $\delta$.

The definition of $\Phi_{X,\spinc}$ depends on a choice of admissible cut
for $X$, but it is not difficult to verify~\cite{HolDiskFour} that
$\Phi_{X,\spinc}$ is independent of this choice, giving a well-defined
four-manifold invariant.

The element $\Phi_{X,\spinc}$ is non-trivial for only finitely many
$\SpinC$ structures over $X$. It vanishes for connected sums of
four-manifolds with $b_2^+(X)>0$,
c.f. Theorem~\ref{HolDiskFour:intro:VanishingTheorem}
of~\cite{HolDiskFour} (compare~\cite{DonaldsonPolynomials} and
\cite{Witten} for corresponding results for Donaldson polynomials and
Seiberg-Witten invariants respectively). In fact, according
to~\cite{HolDiskSymp}, if $(X,\omega)$ is a symplectic four-manifold
$\Phi_{X,\spinccanf}\neq 0$ for the canonical $\SpinC$ structure $\spinccanf$
associated to the symplectic structure.  This
can be seen as an analogue of a theorem of Taubes~\cite{TaubesSympI}
in the Seiberg-Witten context. Whereas Taubes' theorem is proved by
perturbing the Seiberg-Witten equations using a symplectic two-form,
the non-vanishing theorem of $\Phi$ is proved by first associating to
$(X,\omega)$ a compatible Lefschetz pencil, which can be done according to
a theorem of Donaldson, c.f.~\cite{DonaldsonLefschetz}, blowing
up to obtain a Lefschetz fibration, and then analyzing maps between
Floer homology induced by two-handles coming from the singularities in
the Lefschetz fibration, with the help of Theorem~\ref{thm:ExactSeq}.

\section{Knots in $S^3$}
\label{sec:Knots}

We describe here constructions of Heegaard Floer homology applicable
to knots. For simplicity, we restrict attention to knots in
$S^3$. This knot invariant was introduced in~\cite{HolDiskKnots} and
also independently by Rasmussen in~\cite{Rasmussen},
\cite{RasmussenThesis}.  The calculations in
Subsection~\ref{subsec:Calc} are based on the results of
~\cite{AltKnots}, \cite{NoteLens}, \cite{calcKT}. In
Subsection~\ref{subsec:Genus}, we discuss the fact that knot Floer
homology detects the Seifert genus of a knot. This result is proved
in~\cite{GenusBounds}. The relationship with the four-ball genus is
discussed in~\ref{subsec:FourBallGenus}, where we discuss the
concordance invariant of~\cite{FourBall} and~\cite{RasmussenThesis},
and also the method of Owens and Strle~\cite{OwensStrle}. Finally, in
Subsection~\ref{subsec:UnknotOne}, we discuss an application to the
problem of knots with unknotting number one from~\cite{UnknotOne}.
This application uses the Heegaard Floer homology of the branched
double-cover associated to a knot.

\subsection{Knot Floer homology}
\label{subsec:HFK}

In Subsection~\ref{subsec:Knots}, we described a
construction which associates to an oriented knot in a three-manifold
$Y$ a filtration of the chain complex $\CFa(Y)$. Our aim here is to
describe properties of this invariant when the ambient three-manifold
is $S^3$ (although we will be forced to generalize to the case of
knots in a connected sum of copies of $S^2\times S^1$, as we shall see
later). In this case, a knot $K\subset S^3$ induces a filtration of
the chain complex $\CFa(S^3)$, whose homology is a single $\Z$. With
some loss of information, we can take the homology of the associated
graded object, to obtain the ``knot Floer homology''
$$\HFKa_*(K,i)=H_*(\Filt(K,i)/\Filt(K,i-1)).$$ Note that this can be
viewed as one bigraded Abelian group
$$\HFKa(K)=\bigoplus_{d,i\in\Z}\HFKa_d(K,i).$$ We call here $i$ the
filtration level and $d$ (the grading induced from
the Heegaard Floer complex $\CFa(S^3)$) the Maslov grading.

These homology groups satisfy a number of basic properties, which we
outline presently. Sometimes, it is simplest to state these properties 
for $\HFKa(K,i,\Q)$, the homology with rational coefficients:
$\HFKa(K,i,\Q)\cong \HFKa(K,i)\otimes_{\Z}\Q$.

The Euler characteristic is related to the Alexander
polynomial of $K$, $\Delta_K(T)$ by the following formula:
\begin{equation}
\label{eq:EulerChar}
\sum \chi(\HFKa(K,i,\Q))\cm T^{i} =\Delta_K(T)
\end{equation}
(it is interesting to compare this with~\cite{Casson}, \cite{MengTaubes},
and~\cite{FSknots}).
The sign conventions on the Euler characteristic here are given by
\begin{eqnarray*}
\chi(\HFKa(K,i,\Q))
&=& \sum_{d=-\infty}^{+\infty}
(-1)^{d}\cm \rk \left(\HFKa_d(K,i,\Q)\right).
\end{eqnarray*}

Unlike the Alexander polynomial, the knot Floer homology is sensitive
to the chirality of the knot. Specifically,
if
${\overline K}$ denotes the mirror of $K$ (i.e. switch over- and
under-crossings in a projection for $K$), then
\begin{equation}
\label{eq:Reflection}
\HFKa_d(K,i,\Q)\cong \HFKa_{-d}({\overline K},-i,\Q).
\end{equation}
Another symmetry these invariants enjoy is the following conjugation symmetry:
\begin{equation}
\label{eq:KnotConjInv}
\HFKa_d(K,i,\Q)\cong \HFKa_{d-2i}(K,-i,\Q),
\end{equation}
refining the symmetry of the Alexander polynomial.

These groups also satisfy a K\"unneth principle for connected
sums. Specifically, let $K_1$ and $K_2$ be a pair of knots, and let
$K_1\# K_2$
denote their connected sum. Then,
\begin{equation}
\label{eq:KunnethFormula}
\HFKa(K_1\#K_2,i,\Q)\cong
\bigoplus_{i_1+i_2=i} \HFKa(K_1,i_1,\Q)\otimes_\Q
\HFKa(K_2,i_2,\Q)\end{equation}
(see Corollary~\ref{Knots:cor:Kunneth} of~\cite{HolDiskKnots},
and~\cite{RasmussenThesis}).  Of
course, this can be seen as a refinement of the fact that the
Alexander polynomial is multiplicative under connected sums
of knots.

These invariants also satisfy a ``skein exact sequence''
(compare~\cite{FloerTriangles}, \cite{BraamDonaldson}, \cite{AbsGraded},
\cite{Khovanov}). To state it, 
we must generalize to the case of oriented links in $S^3$. This can be
done in the following  manner: an $n$-component oriented link
in $S^3$ gives rise, in a natural way, to an $n$-component oriented
knot in $\#^{n-1}(S^2\times S^1)$. Specifically,
we attach $n-1$ one-handles to
$S^3$, so that the two feet of each one-handle lie on different
components of the link, and so that each link component meets at least
one foot. Next, we form the connected sum of the various components of
the link via standard strips which pass through the one-handles.  In
this way, we view the link invariant for an $n$-component link
$L\subset S^3$ as a knot invariant for the associated knot in
$\#^{n-1}(S^2\times S^1)$.

In this manner, the homology of the associated graded object -- the
link Floer homology -- is a sequence of graded Abelian
groups $\HFKa_*(L,i)$, where here $i\in \Z$. If $L$ has an odd number of
components, the Maslov grading is a $\Z$-grading, while if it has an even
number of components, the Maslov grading takes values in
$\OneHalf+\Z$. As a justification for this convention,
observe that the reflection formula, Equation~\eqref{eq:Reflection},
remains true in the context of links.

\begin{figure}
\mbox{\vbox{\epsfbox{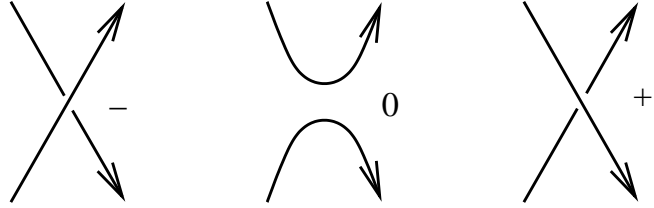}}}
\caption{\label{fig:Skein}
Skein moves at a double-point.}
\end{figure}

Suppose that $L$ is a link, and suppose that $p$ is
a positive crossing of some projection of $L$. Following the usual
conventions from skein theory, there are two other associated links,
$L_0$ and $L_-$, where here $L_-$ agrees with $L_+$, except that the
crossing at $p$ is changed, while $L_0$ agrees with $L_+$, except that
here the crossing $p$ is resolved in a manner consistent
with orientations, as  illustrated
in Figure~\ref{fig:Skein}. There are two cases of the skein exact
sequence, according to whether or not the two strands of $L_+$ which
project to $p$ belong to the same component of $L_+$.

Suppose first that the two strands which project to $p$ belong to the
same component of $L_+$. In this case, the skein exact sequence reads:
\begin{equation}
\label{eq:SkeinExact}
\begin{CD}
...@>>>\HFKa(L_{-})@>>>\HFKa(L_0)@>>>\HFKa(L_{+})@>>>...,
\end{CD}
\end{equation}
where all the maps above respect the splitting of $\HFKa(L)$ into
summands (e.g. $\HFKa(L_-,i)$ is mapped to $\HFKa(L_0,i)$).
Furthermore, the maps to and from $\HFKa(L_0)$ drop degree by
$\OneHalf$. The remaining map  from
$\HFKa(L_+)$ to $\HFKa(L_-)$ does not necessarily respect the 
absolute grading; however, it can be expressed as a sum of
homogeneous maps,
none of which increases absolute grading.
When the two strands belong to
different components, we obtain the following:
\begin{equation}
\label{eq:SkeinExact2}
\begin{CD}
...@>>>\HFKa(L_{-})@>>>\HFKa(L_0)\otimes V @>>>\HFKa(L_{+})@>>>...,
\end{CD}
\end{equation}
where $V$ denotes the four-dimensional vector space $$V=V_{-1}\oplus
V_0\oplus V_1,$$ where here $V_{\pm 1}$ are one-dimensional pieces
supported in degree $\pm 1$, while $V_0$ is a two-dimensional piece
supported in degree $0$.  Moreover, the maps respect the decomposition
into summands, where the $i^{th}$ summand of the middle piece
$\HFKa(L_0)\otimes V$ is given by 
$$\left(\HFKa(L_0,i-1)\otimes V_1\right) 
\oplus \left(\HFKa(L_0,i)\otimes V_0\right) 
\oplus \left(\HFKa(L_0,i+1)\otimes V_{-1}\right).
$$
The shifts in the 
absolute gradings work just as they did in the previous case.

The skein exact sequence is, of course, very closely related to
Theorem~\ref{thm:ExactSeq}. Indeed, its proof proceeds by consider
the surgery long exact sequence associated to an unknot which links
the crossing one is considering, and analyzing the behaviour of the induced
maps, c.f. Section~\ref{Knots:sec:Sequences} of~\cite{HolDiskKnots}.

\subsection{Calculations of knot Floer homology}
\label{subsec:Calc}

It is useful to have a concrete description of the generators
of the knot Floer complex in terms of the combinatorics of a 
knot projection. In fact, the data we fix at first is an
oriented knot projection (with at most double-point singularities), 
equipped with a choice of distinguished edge $e$
which appears in the closure of the unbounded region $A$ in 
the planar projection. We call this data a {\em decorated projection}
for $K$. We denote the planar graph of the projection by $G$.

We can construct a doubly-pointed Heegaard diagram compatible with $K$
from a decorated projection of $K$, as follows.

\begin{figure}
\mbox{\vbox{\epsfbox{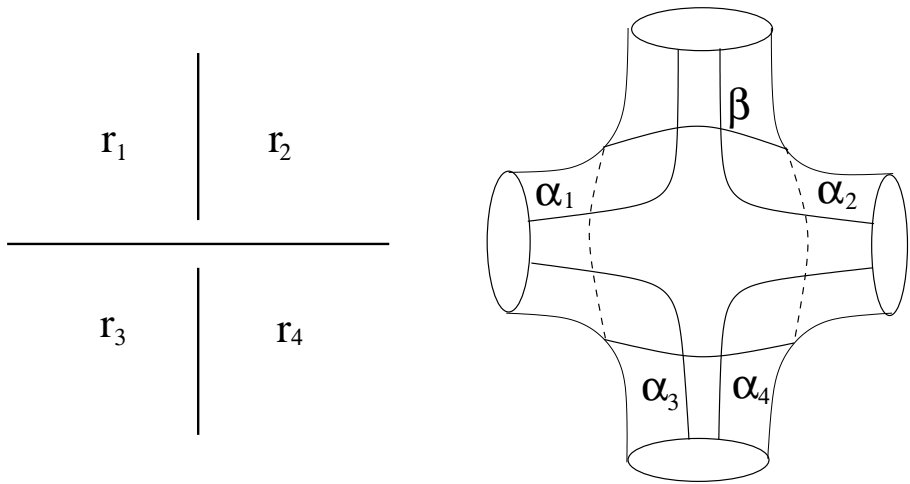}}}
\caption{\label{fig:SpecialHeegaard}
{\bf{Special Heegaard diagram for knot crossings.}} At each crossing
as pictured on the left, we construct a piece of the Heegaard surface
on the right (which is topologically a four-punctured sphere). The
curve $\beta$ is the one corresponding to the crossing on the left;
the four arcs $\alpha_1,...,\alpha_4$ will close up. (Note that if one
of the four regions $r_1,...,r_4$ contains the distinguished edge $e$,
its corresponding $\alpha$-curve should {\em not} be included).
Note that the Heegaard surface is oriented from the outside.}
\end{figure}

Let $B$ denote the other region which contains the edge $e$ in its closure, 
and let
$\Sigma$ be the boundary of a regular neighborhood of $G$,
thought of as a one-complex in $S^3$
(i.e. if our projection has $n$ double-points, then
$\Sigma$ has  genus $n+1$); we orient $\Sigma$ as $\partial
(S^3-\nbd{G})$. We associate to each region $r\in R(G)-A$, an
attaching circle $\alpha_r$ (which follows along the boundary of
$r$). To each crossing $v$ in $G$ we associate an attaching circle
$\beta_v$ as indicated in Figure~\ref{fig:SpecialHeegaard}. In
addition, we let $\mu$ denote the meridian of the knot, chosen to be
supported in a neighborhood of the distinguished edge $e$.

Each vertex $v$ is contained in four (not necessarily distinct)
regions. Indeed, it is clear from Figure~\ref{fig:SpecialHeegaard},
that in a neighborhood of each vertex $v$, there are at most four
intersection points of $\beta_v$ with circles corresponding to the 
four regions which contain $v$. 
(There are fewer than four intersection points with
$\beta_v$ if $v$ is a corner for the unbounded region $A$.) Moreover,
the circle corresponding to $\mu$ meets the circle $\alpha_B$ in a
single point (and is disjoint from the other circles).  Placing
one reference point $w$ and $z$ on each side of $\mu$, we obtain
a doubly-pointed Heegaard diagram for $S^3$ compatible with $K$.

We can now describe the generators $\Ta\cap \Tb$ for the knot Floer
homology in terms of the planar graph $G$ of the projection. 

\begin{defn}
A {\em Kauffman state} (c.f.~\cite{Kauffman}) for a decorated knot
projection of $K$ is a map which associates to each vertex of $G$ one
of the four in-coming quadrants, so that:
\begin{itemize}
\item the quadrants associated to distinct vertices are subsets of distinct
regions in $S^2-G$
\item none of the quadrants is a corner of the distinguished regions
$A$ or $B$ (whose closure contains the edge $e$).
\end{itemize}
\end{defn}

If $K$ is a knot with a decorated projection, it is straightforward to
see that the intersection points $\Ta\cap\Tb$ for the corresponding
Heegaard diagram correspond to Kauffman states for the projection.
Note that Kauffman states have an alternative interpretation, as
maximal trees in the ``black graph'' associated
to a checkerboard coloring of the complement of $G$, c.f.
\cite{Kauffman}.

We can also describe the filtration level and the Maslov grading of a
Kauffman state in combinatorial terms of the decorated knot
projection..

\begin{figure}
\mbox{\vbox{\epsfbox{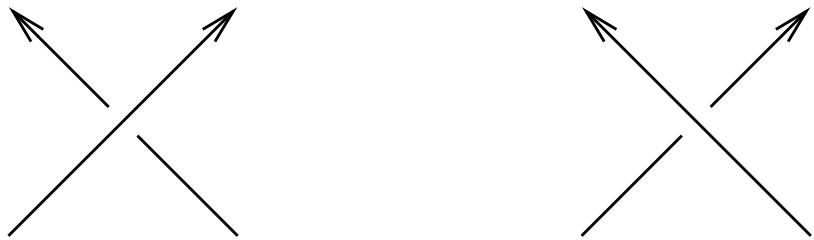}}}
\caption{\label{fig:CrossingSigns}
{\bf{Crossing conventions.}} Crossings of the first kind are assigned
$+1$, and those of the second kind are assigned $-1$.}
\end{figure}

\begin{figure}
\mbox{\vbox{\epsfbox{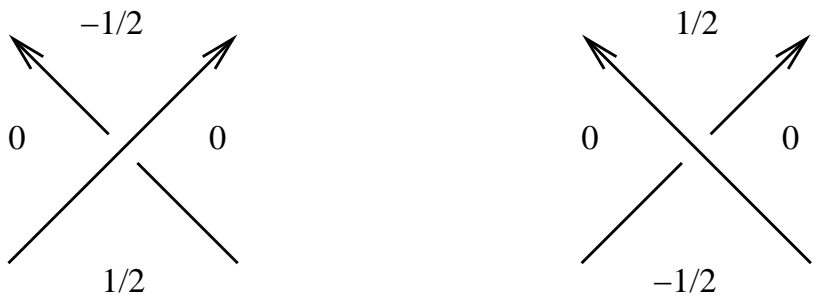}}}
\caption{\label{fig:PointTowards}
{\bf{Local filtration level contributions $s(x,v)$.}} We have
illustrated the local contributions of $s(x,v)$ for both kinds of
crossings. (In both pictures, ``upwards'' region is the one which the
two edges point towards.).}
\end{figure}

To describe the filtration level, note that the orientation on the
knot $K$ associates to each vertex $v\in G$ a distinguished quadrant
whose boundary contains both edges which point towards the vertex
$v$. We call this the quadrant which is ``pointed towards'' at $v$.
There is also a diagonally opposite region which is ``pointed away
from'' (i.e. its boundary contains the two edges pointing away from
$v$).  We define the local filtration contribution of $x$ at $v$, denoted $s(x,v)$, by
the following rule (illustrated in Figure~\ref{fig:PointTowards}),
where $\epsilon(v)$ denotes the sign of the crossing (which we recall
in Figure~\ref{fig:CrossingSigns}): $$2 {\epsilon(v)} s(x,v)=
\left\{
\begin{array}{rl}
1 & {\text{$x(v)$ is the quadrant pointed towards at $v$}}\\ -1 &
{\text{$x(v)$ is the quadrant away from at $v$}}
\\ 0 & {\text{otherwise.}}
\end{array}
\right.$$
The filtration level associated to a Kauffman state, then, is given by the sum
$$\Filt(x)=\sum_{v\in\Vertices(G)}s(x,v).$$ Note that the function $\Filt(x)$
is the $T$-power appearing for the contribution of $x$ to the
symmetrized Alexander polynomial, see~\cite{Alexander},
\cite{KauffmanTwo}.

\begin{figure}
\mbox{\vbox{\epsfbox{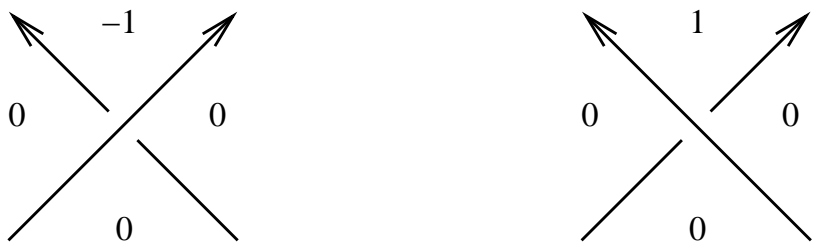}}}
\caption{\label{fig:Grading}
{\bf{Local grading contributions $m(x,v)$.}} We have
illustrated the local contribution of $m(x,v)$.}
\end{figure}

The Maslov grading $\gr(x)$ is defined analogously.  First, at each vertex $v$,
we define the local grading contribution $m(x,v)$.  This local
contributions is non-zero on only one of the four quadrants -- the
one which is pointed away from at $v$.  At this quadrant, the grading
contribution is minus the sign $\epsilon(v)$ of the crossing, as
illustrated in Figure~\ref{fig:Grading}.  Now, the
grading $\gr(x)$ of a Kauffman state $x$ is defined by the formula
$$\gr(x)=\sum_{v\in \Vertices(G)}m(x,v).$$

A verification of these formulas can be found in
Theorem~\ref{AltKnots:thm:States} of~\cite{AltKnots}.

It is clear from the above formulas that if $K$ has an alternating
projection, then $\Filt(\x)-\gr(\x)$ is independent of the choice
of state $\x$. It follows that if we use the chain complex associated
to this Heegaard diagram, then there are no differentials in the knot
Floer homology, and indeed, its rank is determined by its Euler
characteristic. Indeed, by calculating the constant, we get the
following result, proved in Theorem~\ref{AltKnots:thm:KnotHomology}
of~\cite{AltKnots}:

\begin{theorem}
\label{thm:AltKnots}
Let $K\subset S^3$ be an alternating knot in the three-sphere, and
write its symmetrized Alexander polynomial as $$\Delta_K(T) =
a_0+\sum_{s>0} a_s (T^s+T^{-s}),$$ and let $\sigma(K)$ denote its
signature. Then, $\HFKa(S^3,K,s)$ is supported entirely in dimension
$s+\frac{\sigma(K)}{2}$, and indeed $$\HFKa(S^3,K,s)\cong \Z^{|a_s|}.$$
\end{theorem}

Thus, for alternating knots, this choice of Heegaard diagram is
remarkably successful. However, in general, there are differentials
one must grapple with, and these admit, at present, no combinatorial
description in terms of Kauffman states. However, they do respect
certain additional filtrations which can be described in terms of
states, and this property, together with some additional tricks, can
be used to give calculations of knot Floer homology groups in
certain cases, c.f.~\cite{calcKT}, \cite{Eftekhary}. As a particular
example, these filtrations are used in~\cite{calcKT} to show that knot
Floer homology of the eleven-crossing Kinoshita-Terasaka knot (a knot
whose Alexander polynomial is trivial) differs from that of its Conway
mutant.

In a different direction, some knots admit Heegaard diagrams on a
genus one surface.  For these knots, calculation of the differentials
becomes a purely combinatorial matter,
c.f. Section~\ref{Knots:sec:Examples} of~\cite{HolDiskKnots} and
also~\cite{Rasmussen}, \cite{RasmussenThesis}, \cite{GodaMatsudaMorifuji}.

Sometimes, it is more convenient to use more abstract methods to
calculate knot Floer homology. In particular, there is a relationship
between knot Floer homology and the Heegaard Floer homology
three-manifolds obtained by surgery along $K$, c.f.~\cite{HolDiskKnots},
\cite{RasmussenThesis}. With the help of this relationship, we obtain the
following structure for the knot Floer homology of a knot for which
some positive surgery is an $L$-space (proved in Theorem~\ref{NoteLens:thm:FloerHomology} of~\cite{NoteLens}):

\begin{theorem}
\label{thm:LSpaceKnots}
Suppose that $K\subset S^3$ is a knot for which there is a positive
integer $p$ for which $S^3_p(K)$ is an $L$-space. Then, there is an
increasing sequence of non-negative integers
$$n_{-m}<...<n_m$$
with the property that $n_i=-n_{-i}$, with the following significance.
If we let
$$\delta_{i}=\left\{\begin{array}{ll}
0 & {\text{if $i=m$}} \\
\delta_{i+1}-2(n_{i+1}-n_{i})+1 &{\text{if $m-i$ is odd}} \\
\delta_{i+1}-1 & {\text{if $m-i>0$ is even,}}
\end{array}\right.$$
then
$\HFKa(K,j)=0$ unless $j=n_i$ for some $i$, in which case
$\HFKa(K,j)\cong \Z$ and it is supported entirely in dimension $\delta_i$.
\end{theorem}

\begin{figure}
\mbox{\vbox{\epsfbox{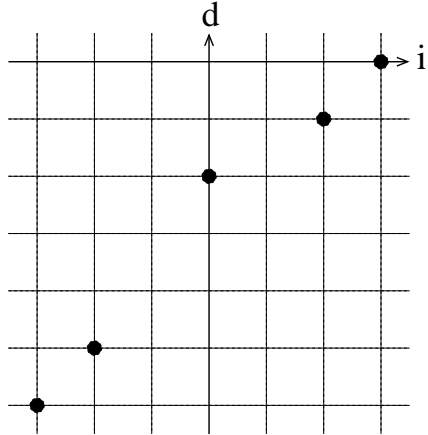}}}
\caption{\label{fig:T34}
Knot Floer homology for the $(3,4)$ torus knot. The dots represent $\Z$ summands, and the bigrading is specified by the $d$ and $i$ coordinates.}
\end{figure}

For example, all (right-handed) torus knots satisfy the hypothesis of
this theorem. (Recal that if $T_{p,q}$ denotes the right-handed
$(p,q)$ torus knot, then $S^3_{pq\pm 1}(T_{p,q})$ is a lens space.)
The knot Floer homology of the $(3,4)$ torus knot is illustrated in
Figure~\ref{fig:T34}. The above theorem can be fruitfully thought of
from three perspectives: as a source of examples of knot Floer
homology calculations (for example, a calculation of the knot Floer
homology of torus knots), as a restriction on knots which admit
$L$-space surgeries (for example, it shows that if $K\subset S^3$
admits a lens space surgery, then all the coefficients of its
Alexander polynomial satisfy $|a_i|\leq 1$), or as a restriction on
$L$-spaces which can arise as surgeries on knots in $S^3$,
c.f.~\cite{NoteLens}.

\subsection{Knot Floer homology and the Seifert genus}
\label{subsec:Genus}

A knot $K\subset S^3$ can be realized as the boundary of
an embedded, orientable surface in $S^3$. Such a surface is called a
Seifert surface for $K$, and the minimal genus of any Seifert surface
for $K$ is called its {\em Seifert genus}, denoted $g(K)$. Of course,
a knot has $g(K)=0$ if and only if it is the unknot.

The knot Floer homology of $K$ detects the Seifert genus, and in
particular it distinguishes the unknot, according to the following
result proved in~\cite{GenusBounds}. To state it, we first define the
degree of the knot Floer homology to be the integer
$$\deg\HFKa(K) = \max\{i\in\Z\big| \HFKa(K,i)\neq 0\}.$$

\begin{theorem}
\label{thm:GenusBounds}
For any knot $K\subset S^3$, $g(K)=\deg\HFKa(K)$.
\end{theorem}

Given a Seifert surface of genus $g$ for $K$, one can construct a
Heegaard diagram for which all the points in $\Ta\cap \Tb$ have
filtration level $\leq g$. This gives at once the bound
$$\deg\HFKa(K)\leq g(K)$$ (this result is analogous
to a classical bound on the genus of a knot in terms of the degree of
its Alexander polynomial).

The inequality in the other direction is much more subtle, involving
much of the theory described so far. First, one relates the degree of
the knot Floer homology by a similar quantity defined using the Floer
homology of the zero-surgery $S^3_0(K)$. Next, one appeals to a
theorem of Gabai~\cite{GabaiKnots}, according to which if $K$ is a
knot with Seifert genus $g>0$, then $S^3_{0}(K)$ admits a taut foliation
${\mathcal F}$ whose first Chern class is $g-1$ times a generator for
$H^2(S^3_0(K);\Z)$.  The taut foliation naturally induces a symplectic
structure on $[-1,1]\times S^3_0(K)$, according to a result of
Eliashberg and Thurston~\cite{EliashbergThurston}, which, according to
a recent result of Eliashberg~\cite{Eliashberg}, \cite{Etnyre} can be
embedded in a closed symplectic four-manifold $X$ (indeed, one can
arrange for $S^3_0(K)$ to divide the four-manifold $X$ into two pieces
with $b_2^+(X_i)>0$). The non-vanishing of the four-manifold invariant
$\Phi_{X,\spinccanf}$ for a symplectic four-manifold can then be used to
prove that the Heegaard Floer homology of $S^3_0(K)$ is non-trivial in
the $\SpinC$ structure gotten by restricting the canonical $\SpinC$
structure $\spinccanf$ of the ambient symplectic
four-manifold -- i.e. this is the $\SpinC$
structure belonging to the foliation ${\mathcal F}$. The details of
this argument are given in~\cite{GenusBounds}.

Since the generators of knot Floer homology can be thought of from a
Morse-theoretic point of view as simultaneous trajectories of gradient
flow-lines, Theorem~\ref{thm:GenusBounds} immediately gives a curious
Morse-theoretic characterization of the Seifert genus of $K$, as
the minimum over all Morse functions compatible with $K$ of the
maximal filtration level of any simultaneous trajectory.

\subsection{The four-ball genus}
\label{subsec:FourBallGenus}

A knot $K\subset S^3$ can be viewed as a knot in the boundary of the
four-ball, and as such, it can be realized as the boundary of a smoothly
embedded oriented
surface in the four-ball. The minimal genus of any such
surface is called the {\em four-ball genus} of the knot, and it is 
denoted $g^*(K)$. Obviously, $g^*(K)\leq g(K)$. In  general, 
$g^*(K)$ is
quite difficult to calculate.

Lower bounds on $g^*(K)$ can be obtained from Heegaard Floer homology.
The construction involves going deeper into the knot filtration.
Specifically, as explained explained in Subsection~\ref{subsec:Knots},  the
filtered chain homotopy type of the sequence of inclusions
$$
.....\subseteq \Filt(K,i)\subseteq
\Filt(K,i+1)\subseteq...\subseteq\CFa(S^3)
$$
is a knot invariant; passing to the homology of the associated graded
object constitutes some loss of information. 
There is a quantity associated to the filtered complex which
goes beyond knot Floer homology, and that is the integer $\tau(K)$ which
is defined by
$$\tau(K)=\min\{i\in\Z\big|
H_*(\Filt(K,i))\longrightarrow \HFa(S^3)
~\text{is non-trivial}\}.$$

It is proved in~\cite{FourBall}, \cite{RasmussenThesis} that
\begin{equation}
\label{eq:FourBallGenus}
|\tau(K)|\leq g^*(K).
\end{equation}

The above inequality can be used to prove a property of $\tau$ which
underscores its analogy with the correction terms $d(Y,\spinc)$
described earlier. To put this result into context, we give a
definition.  Two knots $K_1$ and $K_2$ are said to be {\em concordant}
if there is a smoothly embedded cylinder $C$ in $[1,2]\times S^3$
with $C\cap \{i\}\times S^3 = K_i$ for $i=1,2$. The set of concordance
classes of knots can be made into an Abelian group, under the
connected
sum operation. It follows from Inequality~\eqref{eq:FourBallGenus}, together
with the additivity of $\tau$ under connected sums, that 
$\tau$  gives a homomorphism from the concordance group of knots to $\Z$.

Another such homomorphism is provided by $\sigma(K)/2$. 
In fact, according to Theorem~\ref{thm:AltKnots}, 
\begin{equation}
\label{eq:AltTau}
2\tau(K)=-\sigma(K)
\end{equation}
when $K$ is alternating. By contrast, Theorem~\ref{thm:LSpaceKnots}
gives many examples where Equation~\eqref{eq:AltTau} fails. Indeed,
combining Theorem~\ref{thm:LSpaceKnots} with
Theorem~\ref{thm:GenusBounds} and Equation~\eqref{eq:EulerChar}
and~\eqref{eq:FourBallGenus}, we see that if $K$ is a knot which
admits a positive surgery which is an $L$-space surgery, then $$
\tau(K)=g(K)=g^*(K)=\deg\Delta_K;$$
and in particular, if $K=T_{p,q}$, then
\begin{equation}
\label{eq:TauPQ}
\tau(K)=\frac{(p-1)(q-1)}{2}.
\end{equation}
Note that the fact that $g^*(T_{p,q})$ is given by the above formula
was conjectured by Milnor and first proved by
Kronheimer and Mrowka using gauge theory~\cite{KMMilnor}.

In~\cite{LivingstonCalc}, see also~\cite{LivingstonConc} Livingston shows that properties of the
concordance invariant $\tau(K)$ (specifically, the fact that it is a
homomorphism whose absolute value bounds the four-ball genus of $K$,
and satisfies Equation~\eqref{eq:TauPQ}) leads to the result that if
$K$ is the closure of a positive on $k$ strands with $n$ crossings,
then $$\tau(K)=\frac{n-k+1}{2}=g^*(K).$$ The second of these equations
was proved first by Rudolph~\cite{Rudolph} using the local
Thom conjecture proved by Kronheimer and Mrowka~\cite{KMMilnor}.
Further links between the Thurston-Bennequin invariant and $\tau$ are
explored by Plamenevskaya, see~\cite{PlamenevskayaTB}.

A different method for bounding $g^*(K)$ is given by Owens and Strle
in~\cite{OwensStrle}, where they describe a method using the
correction terms for the branched double-cover of $S^3$ along $K$,
$\Sigma(K)$. Under favorable circumstances, their method gives an
obstruction for Murasugi's bound 
\begin{equation}
\label{eq:MurasugiBound}
|\sigma(K)|\leq 2g^*(K).
\end{equation}
to being sharp. Specifically, taking the branched double-cover
$\Sigma(F)$ of a surface $F$ in $B^4$ which bounds $K$, one obtains a
four-manifold which bounds $\Sigma(K)$.  When $F$ is a surface with
$2g(F)=\sigma(K)$, the branched double $\Sigma(F)$ is a four-manifold
with definite intersection form, whose second Betti number is
$2g(F)$. Comparing this construction with the inequality from
Theorem~\ref{thm:IntForms}, one obtains restrictions on the correction
terms of $\Sigma(K)$, which sometimes can rule out the existence
of such surfaces $F$.

\subsection{Unknotting numbers}
\label{subsec:UnknotOne}

Recall that the {\em unknotting number} of a knot $K$, denoted $u(K)$,
is the the minimal number of crossing-changes required to unknot
$K$. An unknotting of $K$ can be realized as an immersed disk in $B^4$
which bounds $S^3$. Resolving the self-intersections, one gets the
inequality $g^*(K) \leq u(K)$.  However, there are circumstances where
one needs better bounds (most strikingly, for any non-trivial slice
knot). In~\cite{UnknotOne}, we describe an obstruction to knots $K$
having $u(K)=1$, in terms of the correction terms of the branched
double-cover of $K$.

This construction works best in the case where $K$ is alternating.  In
this case, the branched double-cover $\Sigma(K)$ is an $L$-space,
c.f.~\cite{BrDCov}.  A classical construction of
Montesinos~\cite{Montesinos} shows that if $u(K)=1$, then $\Sigma(K)$
can be obtained as $\pm D/2$ surgery on another knot $C$ in $S^3$,
where here $D$ denotes the determinant of the knot $K$
($D=|\Delta_{K}(-1)|$). On the one hand, correction terms for an
$L$-space which is realized as $n/2$ surgery (for some integer $n$) on
a knot in $S^3$ satisfy certain symmetries
(c.f. Theorem~\ref{UnknotOne:thm:LSpaceSymmetry} of~\cite{UnknotOne});
on the other hand, the correction terms of the branched double-cover
of an alternating knot can be calculated explicitly by classical data
associated to an alternating projection of $K$, c.f.~\cite{BrDCov}.
Rather than recalling the result here, we content ourselves with
illustrating an alternating knot $K$ (listed as $8_{10}$ in the
Alexander-Briggs notation) whose unknotting number was previously unknown,
but which can now be shown to have $u(K)=2$ using
these techniques.

\begin{figure}
\mbox{\vbox{\epsfbox{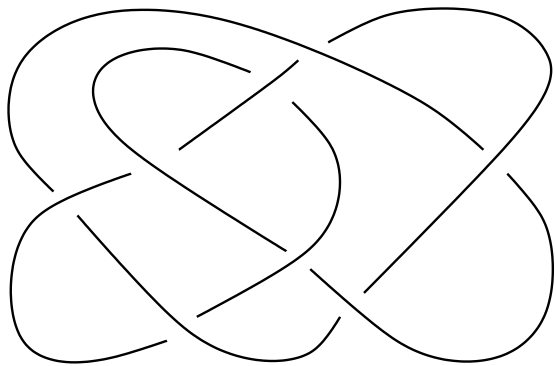}}}
\caption{\label{fig:8s10}
{\bf A knot with $u(K)=2$}.}
\end{figure}

\section{Problems and Questions}
\label{sec:Questions}

The investigation of Heegaard Floer homology naturally leads us to 
the following problems and questions.

Perhaps the most important problem in this circle of ideas is the
following:

\vskip.2cm
\noindent{\bf{Problem 1:}} {\em Give a purely
combinatorial calculation of  Heegaard Floer homology or,
more generally, the Heegaard Floer functor.}
\vskip.2cm

In certain special cases, combinatorial calculations can be given, for
example~\cite{SomePlumbs}, \cite{AltKnots}. This problem would be very
interesting to solve even for certain restricted classes of
three-manifolds, for example for those which fiber over the circle,
compare.~\cite{HutchingsSullivan}, \cite{SeidelDehn}, \cite{SeidelMutation}
\cite{DonaldsonSmith}.

In a related direction, it is simpler to consider the case
of knots in $S^3$. Recall that in Section~\ref{sec:Knots}
we showed that the generators of the knot Floer complex can be
thought of as Kauffman states.

\vskip.2cm
\noindent{\bf{Question 2:}} {\em Is there a combinatorial description
of the differential on Kauffman states whose homology gives the knot
Floer homology of $K$.}
\vskip.2cm

It is intriguing to compare this with Khovanov's new invariants for
links, see~\cite{Khovanov}. These invariants have a very similar
structure to the knot Floer homology considered here, except that
their Euler characteristic gives the Jones polynomial, see
also~\cite{Jacobsson},
\cite{TangleCobordism}, \cite{KhovanovSL3}, \cite{KhovanovRozansky}
\cite{BarNatan}. Indeed, the similarities are further underscored by
the work of E-S. Lee~\cite{EunSooLee}, who describes a spectral
sequence which converges to a vector space of dimension $2^n$, where
$n$ is the number of components of $L$, see also~\cite{RasmussenKhovanov}.

In a similar vein, it would be interesting to give combinatorial
calculations of the numerical invariants arising from Heegaard Floer
homology, specifically, the correction terms $d(Y,\spinc)$ or the
concordance invariant $\tau(K)$. An intriguing conjecture of
Rasmussen~\cite{RasmussenKhovanov} relates $\tau$ with a numerical
invariant coming from Khovanov homology.

\vskip.2cm
\noindent{\bf Problem 3:} {\em Establish the conjectured relationship
between Heegaard Floer homology and Seiberg-Witten theory.}
\vskip.2cm

There are two approaches one might take to this problem. One direct,
analytical approach would be to analyze moduli spaces of solutions to
the Seiberg-Witten equations over a three-manifold
equipped with a Heegaard decomposition. Another approach
would involve an affirmative answer to the following question:

\vskip.2cm
\noindent{\bf  Question 4:} {\em Is there an axiomatic characterization
of Heegaard Floer homology?}
\vskip.2cm

A {\em Floer functor} is a map which associates to any closed,
oriented three-manifold $Y$ a $\Zmod{2}$-graded Abelian group
${\mathcal H}(Y)$ and to any cobordism $W$ from $Y_1$ to $Y_2$ a
homorphism ${\mathcal D}_W\colon {\mathcal H}(Y_1) \longrightarrow
{\mathcal H}(Y_2)$, which is natural under composition of cobordisms,
and which induce exact sequences for triples of three-manifolds $(Y_0,
Y_1, Y_2)$ which are related as in the hypotheses of
Theorem~\ref{thm:ExactSeq}. It is interesting to observe that if ${\mathcal T}$
is a natural transformation between Floer functors ${\mathcal H}$ and
${\mathcal D}$ to ${\mathcal H}'$ and ${\mathcal D}'$, then if ${\mathcal T}$
induces an isomorphism ${\mathcal T}(S^3) \colon {\mathcal H}(S^3)
\longrightarrow {\mathcal H}'(S^3)$, then ${\mathcal T}$ induces isomorphisms
for all three-manifolds ${\mathcal T}(Y) \colon {\mathcal H}(Y)
\longrightarrow {\mathcal H}'(Y)$. This can be proved from Kirby
calculus, following the outline laid out in~\cite{BraamDonaldson}. We
know from~\cite{KMOS} that monopole Floer homology (taken with
$\Zmod{2}$ coefficients) is a Floer functor in this sense, and also
(Theorem~\ref{thm:ExactSeq}) that Heegaard Floer homology is a Floer
functor. Unfortunately, this still falls short of giving an axiomatic
characterization: one needs axioms which are sufficient to assemble a
natural transformation ${\mathcal T}$.

\vskip.2cm
\noindent{\bf  Problem 5:} {\em Develop cut-and-paste techniques for calculating 
the Heegaard Floer homology of $Y$ in terms of data associated to its pieces.}
\vskip.2cm

As a special case, one can ask how the knot Floer homology of a
satellite knot can be calculated from data associated to the companion
and the pattern.  Of course, the K\"unneth principle for connected
sums can be viewed as an example of this. Another example, of
Whitehead doubling, has been studied by Eftekhary~\cite{Eftekhary}.

We have seen that the set of $\SpinC$ structures for which the
Heegaard Floer homology of $Y$ is non-trivial determines the Thurston
norm of $Y$. It is natural to ask what additional topological information
is contained in the groups themselves. It is possible that these groups
contain further information about foliations over $Y$. 

We specialize to the case of knot Floer homology. If $K\subset S^3$ is
a fibered knot of genus $g$, then it is shown in~\cite{HolDiskContact}
that $\HFKa(K,g)\cong \Z$.

\vskip.2cm
{\noindent{\bf Question 6:} {\em If $K\subset S^3$ is a knot with
genus $g$ and $\HFKa(K,g)\cong \Z$, does it follow that $K$ is
fibered?}}
\vskip.2cm

Calculations give some evidence that the answer to the above 
question is positive.

\vskip.2cm
\noindent{\bf Question 7:} {\em If $K\subset S^3$ is a knot,
is there an explicit relationship
between the fundamental group of $S^3-K$ and
the knot Floer homology $\HFKa(K)$?}
\vskip.2cm

The above question is, of course, very closely related to the following:

\vskip.2cm
\noindent{\bf Question 8:} {\em Is there an explicit relationship
between the Heegaard Floer homology and the fundamental group of $Y$?}
\vskip.2cm

For example, one could try to relate the Heegaard Floer homology with
the instanton Floer homology of $Y$. (Note, though, that presently
instanton Floer homology is defined only for a restricted class of
three-manifolds, c.f.~\cite{BraamDonaldson}.) A link between
Seiberg-Witten theory and instanton Floer homology is given by
Pidstrygach and Tyurin's $PU(2)$ monopole
equations~\cite{FeehanLeness}. This connection has been successfully
exploited in Kronheimer and Mrowka's recent proof that all knots in
$S^3$ have Property P~\cite{KMpropP}.

The conjectured relationship with Seiberg-Witten invariants raises
further questions.  Specifically, Bauer and Furuta ~\cite{Furuta},
~\cite{BauerFuruta} have constructed refinements of the Seiberg-Witten
invariant which use properties of the Seiberg-Witten equations beyond
merely their solution counts. Correspondingly, these invariants carry
topological information about four-manifolds beyond their usual
Seiberg-Witten invariants.  A three-dimensional analogue is studied in
work of Manolescu and Kronheimer, see~\cite{Manolescu},
\cite{ManolescuKronheimer}

\vskip.2cm
\noindent{\bf Question 9:} {\em  Is there a refinement of 
the four-manifold invariant $\Phi$ defined using Heegaard Floer
homology which captures the information in the Bauer-Furuta construction?}
\vskip.2cm

In the opposite direction, it is natural to study the following:

\vskip.2cm
\noindent{\bf Problem  10:} {\em Construct a gauge-theoretic
analogue of knot Floer homology.}
\vskip.2cm

Recall from Section~\ref{sec:BasProp} that there is a class of
three-manifolds whose Heegaard Floer homology is as simple as
possible, the so-called $L$-spaces.  This class of three-manifolds
includes all lens spaces, and more generally branched double-covers of
alternating knots.  The set of $L$-spaces is closed under connected
sums. According to~\cite{GenusBounds}, $L$-spaces admit no
(coorientable) taut foliations. A striking theorem of
N{\'e}methi~\cite{Nemethi} characterizes those $L$-spaces which are
boundaries of negative-definite plumbings of spheres: they are the
links of rational surface singularities.  Note also that there is an
analogous class of three-manifolds in the context of Seiberg-Witten
monopole Floer homology, c.f.~\cite{KMOS}.

\vskip.2cm
\noindent{\bf Question  11:} {\em Is there a topological
characterization of $L$-spaces (i.e. which makes no reference to 
Floer homology)?}

\commentable{
\bibliographystyle{plain}
\bibliography{biblio}

\begin{thebibliography}{100}

\bibitem{Casson}
S.~Akbulut and J.~D. McCarthy.
\newblock {\em Casson's invariant for oriented homology {$3$}-spheres},
  volume~36 of {\em Mathematical Notes}.
\newblock Princeton University Press, Princeton, NJ, 1990.
\newblock An exposition.

\bibitem{Alexander}
J.~W. Alexander.
\newblock Topological invariants of knots and links.
\newblock {\em Trans. Amer. Math. Soc.}, 30(2):275--306, 1928.

\bibitem{AtiyahFloer}
M.~F. Atiyah.
\newblock {\em Floer homology}, pages 105--108.
\newblock Number 133 in Progr. Math. Birkh{\"a}user, 1995.

\bibitem{APSII}
M.~F. Atiyah, V.~K. Patodi, and I.~M. Singer.
\newblock Spectral asymmetry and {R}iemannian geometry. {II}.
\newblock {\em Math. Proc. Cambridge Philos. Soc}, 78(3):405--432, 1975.

\bibitem{BarNatan}
D.~Bar-Natan.
\newblock On {K}hovanov's categorification of the {J}ones polynomial.
\newblock {\em Algebraic and Geometric Topology}, 2:337--370, 2002.

\bibitem{BauerFuruta}
S.~Bauer and M.~Furuta.
\newblock A stable cohomotopy refinement of {Seiberg-Witten} invariants: {I}.
\newblock math.DG/0204340.

\bibitem{BraamDonaldson}
P.~Braam and S.~K. Donaldson.
\newblock Floer's work on instanton homology, knots, and surgery.
\newblock In H.~Hofer, C.~H. Taubes, A.~Weinstein, and E.~Zehnder, editors,
  {\em The Floer Memorial Volume}, number 133 in Progress in Mathematics, pages
  195--256. Birkh{\"a}user, 1995.

\bibitem{CGLS}
M.~Culler, C.~McA. Gordon, J.~Luecke, and P.~B. Shalen.
\newblock Dehn surgery on knots.
\newblock {\em Ann. of Math}, 125(2):237--300, 1987.

\bibitem{DonaldsonSmith}
S.~Donaldson and I.~Smith.
\newblock Lefschetz pencils and the canonical class for symplectic
  four-manifolds.
\newblock {\em Topology}, 42(4):743--785, 2003.

\bibitem{Donaldson}
S.~K. Donaldson.
\newblock An application of gauge theory to four-dimensional topology.
\newblock {\em J. Differential Geom.}, 18(2):279--315, 1983.

\bibitem{Chambers}
S.~K. Donaldson.
\newblock Irrationality and the {$h$}-cobordism conjecture.
\newblock {\em J. Differential Geom.}, 26(1):141--168, 1987.

\bibitem{DonaldsonPolynomials}
S.~K. Donaldson.
\newblock Polynomial invariants for smooth four-manifolds.
\newblock {\em Topology}, 29(3):257--315, 1990.

\bibitem{DonaldsonSurvey}
S.~K. Donaldson.
\newblock The {S}eiberg-{W}itten equations and {$4$}-manifold topology.
\newblock {\em Bull. Amer. Math. Soc. (N.S.)}, 33(1):45--70, 1996.

\bibitem{DonaldsonLefschetz}
S.~K. Donaldson.
\newblock Lefschetz pencils on symplectic manifolds.
\newblock {\em J. Differential Geom.}, 53(2):205--236, 1999.

\bibitem{DonaldsonBook}
S.~K. Donaldson.
\newblock {\em Floer homology groups in {Y}ang-{M}ills theory}, volume 147 of
  {\em Cambridge Tracts in Mathematics}.
\newblock Cambridge University Press, Cambridge, 2002.
\newblock With the assistance of M. Furuta and D. Kotschick.

\bibitem{DonKron}
S.~K. Donaldson and P.~B. Kronheimer.
\newblock {\em The Geometry of Four-Manifolds}.
\newblock Oxford Mathematical Monographs. Oxford University Press, 1990.

\bibitem{DostoglouSalamon}
S.~Dostoglou and D.~A. Salamon.
\newblock Self-dual instantons and holomorphic curves.
\newblock {\em Ann. of Math. (2)}, 139(3):581--640, 1994.

\bibitem{Eftekhary}
E.~Eftekhary.
\newblock Knot {F}loer homologies for pretzel knots.
\newblock math.GT/0311419.

\bibitem{Eliashberg}
Y.~M. Eliashberg.
\newblock Few remarks about symplectic filling.
\newblock math.SG/0311459.

\bibitem{EliashbergThurston}
Y.~M. Eliashberg and W.~P. Thurston.
\newblock {\em Confoliations}.
\newblock Number~13 in University Lecture Series. American Mathematical
  Society, 1998.

\bibitem{Elkies}
N.~D. Elkies.
\newblock A characterization of the {${Z}\sp n$} lattice.
\newblock {\em Math. Res. Lett.}, 2(3):321--326, 1995.

\bibitem{Etnyre}
J.~B. Etnyre.
\newblock On symplectic fillings.
\newblock math.SG/0312091, 2003.

\bibitem{FeehanLeness}
P.~M.~N. Feehan and T.~G. Leness.
\newblock {$\rm SO(3)$} monopoles, level-one {S}eiberg-{W}itten moduli spaces,
  and {W}itten's conjecture in low degrees.
\newblock {\em Topology Appl.}, 124(2):221--326, 2002.

\bibitem{FSsfs}
R.~Fintushel and R.~J. Stern.
\newblock Seifert fibered homology three-spheres.
\newblock {\em Proc. of the London Math. Soc.}, 61:109--137, 1990.

\bibitem{FSknots}
R.~Fintushel and R.~J. Stern.
\newblock Knots, links, and {$4$}-manifolds.
\newblock {\em Invent. Math.}, 134(2):363--400, 1998.

\bibitem{InstantonFloer}
A.~Floer.
\newblock An instanton-invariant for 3-manifolds.
\newblock {\em Comm. Math. Phys.}, 119:215--240, 1988.

\bibitem{FloerLag}
A.~Floer.
\newblock Morse theory for {L}agrangian intersections.
\newblock {\em J. Differential Geometry}, 28:513--547, 1988.

\bibitem{FloerTriangles}
A.~Floer.
\newblock Instanton homology and {D}ehn surgery.
\newblock In H.~Hofer, C.~H. Taubes, A.~Weinstein, and E.~Zehnder, editors,
  {\em The Floer Memorial Volume}, number 133 in Progress in Mathematics, pages
  77--97. Birkh{\"a}user, 1995.

\bibitem{FloerHofer}
A.~Floer and H.~Hofer.
\newblock Coherent orientations for periodic orbit problems in symplectic
  geometry.
\newblock {\em Math. Z.}, 212(1):13--38, 1993.

\bibitem{FloerHoferSalamon}
A.~Floer, H.~Hofer, and D.~Salamon.
\newblock Transversality in elliptic {M}orse theory for the symplectic action.
\newblock {\em Duke Math. J}, 80(1):251--29, 1995.

\bibitem{FriedmanMorgan}
R.~Friedman and J.~W. Morgan.
\newblock {\em Smooth four-manifolds and complex surfaces}, volume~27 of {\em
  Ergebnisse der Mathematik und ihrer Grenzgebiete (3)}.
\newblock Springer-Verlag, Berlin, 1994.

\bibitem{Froyshov}
K.~A. Fr{\o}yshov.
\newblock The {S}eiberg-{W}itten equations and four-manifolds with boundary.
\newblock {\em Math. Res. Lett}, 3:373--390, 1996.

\bibitem{FOOO}
K.~Fukaya, Y-G. Oh, K.~Ono, and H.~Ohta.
\newblock {\em Lagrangian intersection Floer theory---anomaly and obstruction}.
\newblock Kyoto University, 2000.

\bibitem{Furuta}
M.~Furuta.
\newblock Finite dimensional approximations in geometry.
\newblock In {\em Proceedings of the International Congress of Mathematicians,
  Vol. II (Beijing, 2002)}, pages 395--403, Beijing, 2002. Higher Ed. Press.

\bibitem{GabaiKnots}
D.~Gabai.
\newblock Foliations and the topology of {$3$}-manifolds {III}.
\newblock {\em J. Differential Geom.}, 26(3):479--536, 1987.

\bibitem{GodaMatsudaMorifuji}
H.~Goda, H.~Matsuda, and T.~Morifuji.
\newblock Knot {F}loer homology of {$(1,1)$}-knots.
\newblock math.GT/0311084.

\bibitem{GompfStipsicz}
R.~E. Gompf and A.~I. Stipsicz.
\newblock {\em {$4$}-manifolds and Kirby calculus}, volume~20 of {\em Graduate
  Studies in Mathematics}.
\newblock American Mathematical Society, 1999.

\bibitem{GordonLueckeI}
C.~McA. Gordon and J.~Luecke.
\newblock Knots are determined by their complements.
\newblock {\em J. Amer. Math. Soc.}, 2(2):371--415, 1989.

\bibitem{Gromov}
M.~Gromov.
\newblock Pseudo holomorphic curves in symplectic manifolds.
\newblock {\em Inventiones Mathematicae}, 82:307--347, 1985.

\bibitem{HutchingsSullivan}
M.~Hutchings and M.~Sullivan.
\newblock The periodic {F}loer homology of a {D}ehn twist.
\newblock http://math.berkeley.edu/~hutching/pub/index.htm, 2002.

\bibitem{IonelParker}
E.~Ionel and T.~H. Parker.
\newblock Relative {G}romov-{W}itten invariants.
\newblock {\em Ann. of Math. (2)}, 157(1):45--96, 2003.

\bibitem{Jacobsson}
M.~Jacobsson.
\newblock An invariant of link cobordisms from {K}hovanov's homology theory.
\newblock math.GT/0206303, 2002.

\bibitem{Kauffman}
L.~H. Kauffman.
\newblock {\em Formal knot theory}.
\newblock Number~30 in Mathematical Notes. Princeton University Press, 1983.

\bibitem{KauffmanTwo}
L.~H. Kauffman.
\newblock {\em On knots}.
\newblock Number 115 in Annals of Mathematics Studies. Princeton University
  Press, 1987.

\bibitem{KhovanovSL3}
M.~Khovanov.
\newblock $sl(3)$ link homology {I}.
\newblock math.QA/0304375.

\bibitem{Khovanov}
M.~Khovanov.
\newblock A categorification of the {J}ones polynomial.
\newblock {\em Duke Math. J.}, 101(3):359--426, 2000.

\bibitem{TangleCobordism}
M.~Khovanov.
\newblock An invariant of tangle cobordisms.
\newblock math.QA/0207264, 2002.

\bibitem{KhovanovRozansky}
M.~Khovanov and L.~Rozansky.
\newblock Matrix factorizations and link homology.
\newblock math.QA/0401268.

\bibitem{ManolescuKronheimer}
P.~B. Kronheimer and C.~Manolescu.
\newblock Periodic {F}loer pro-spectra from the {S}eiberg-{W}itten equations.
\newblock math.GT/0203243.

\bibitem{KMbook}
P.~B. Kronheimer and T.~S. Mrowka.
\newblock Floer homology for {S}eiberg-{W}itten {M}onopoles.
\newblock In preparation.

\bibitem{KMMilnor}
P.~B. Kronheimer and T.~S. Mrowka.
\newblock Gauge theory for embedded surfaces. {I}.
\newblock {\em Topology}, 32(4):773--826, 1993.

\bibitem{KMthom}
P.~B. Kronheimer and T.~S. Mrowka.
\newblock The genus of embedded surfaces in the projective plane.
\newblock {\em Math. Research Letters}, 1:797--808, 1994.

\bibitem{KMPolyStruct}
P.~B. Kronheimer and T.~S. Mrowka.
\newblock Embedded surfaces and the structure of {D}onaldson's polynomial
  invariants.
\newblock {\em J. Differential Geometry}, pages 573--734, 1995.

\bibitem{KMcontact}
P.~B. Kronheimer and T.~S. Mrowka.
\newblock Monopoles and contact structures.
\newblock {\em Invent. Math.}, 130(2):209--255, 1997.

\bibitem{KMpropP}
P.~B. Kronheimer and T.~S. Mrowka.
\newblock Witten's conjecture and {P}roperty {P}.
\newblock math.GT/0311489, 2003.

\bibitem{KMOS}
P.~B. Kronheimer, T.~S. Mrowka, P.~S. Ozsv{\'a}th, and Z.~Szab{\'o}.
\newblock Monopoles and lens space surgeries.
\newblock math.GT/0310164.

\bibitem{EunSooLee}
E.~S. Lee.
\newblock The support of the {K}hovanov's invariants for alternating knots.
\newblock math.GT/0201105, 2002.

\bibitem{LiRuan}
A-M. Li and Y.~Ruan.
\newblock Symplectic surgery and {G}romov-{W}itten invariants of {C}alabi-{Y}au
  3-folds.
\newblock {\em Invent. Math.}, 145(1):151--218, 2001.

\bibitem{LiscaStipsicz}
P.~Lisca and A.~I. Stipsicz.
\newblock Heegaard {F}loer invariants and tight contact three-manifolds.
\newblock math.SG/0303280, 2003.

\bibitem{LivingstonCalc}
C.~Livingston.
\newblock Computations of the {Ozsv\'ath-Szab\'o} concordance invariant.
\newblock math.GT/0311036.

\bibitem{LivingstonConc}
C.~Livingston.
\newblock Splitting the concordance group of algebraically slice knots.
\newblock {\em Geom. Topol.}, 7:641--463, 2003.

\bibitem{Manolescu}
C.~Manolescu.
\newblock Seiberg-{W}itten-{F}loer stable homotopy type of three-manifolds with
  $b_1=0$.
\newblock {\em Geom. Topol.}, 7:889--932, 2003.

\bibitem{MarcolliWang}
M.~Marcolli and B-L. Wang.
\newblock Equivariant {S}eiberg-{W}itten {F}loer homology.
\newblock {\em Comm. Anal. Geom.}, 9(3):451--639, 2001.

\bibitem{MengTaubes}
G.~Meng and C.~H. Taubes.
\newblock {\underline{SW}}={M}ilnor torsion.
\newblock {\em Math. Research Letters}, 3:661--674, 1996.

\bibitem{Montesinos}
J.~M. Montesinos.
\newblock Surgery on links and double branched covers of {$S\sp{3}$}.
\newblock In {\em Knots, groups, and $3$-manifolds (Papers dedicated to the
  memory of R. H. Fox)}, pages 227--259. Ann. of Math. Studies, No. 84.
  Princeton Univ. Press, Princeton, N.J., 1975.

\bibitem{Morgan}
J.~W. Morgan.
\newblock {\em The {S}eiberg-{W}itten Equations and Applications to the
  Topology of Smooth Four-Manifold}.
\newblock Number~44 in Mathematical Notes. Princeton University Press, 1996.

\bibitem{MMR}
J.~W. Morgan, T.~S. Mrowka, and D.~Ruberman.
\newblock {\em The {$L^2$}-Moduli Space and a Vanishing Theorem for Donaldson
  Polynomial Invariants}.
\newblock Number~II in Monographs in Geometry and Topology. International
  Press, 1994.

\bibitem{MOY}
T.~S. Mrowka, P.~S. Ozsv{\'a}th, and B.~Yu.
\newblock Seiberg-{W}itten monopoles on {S}eifert fibered spaces.
\newblock {\em Comm. in Analysis and Geometry}, 5(4):685--793, 1997.

\bibitem{Nemethi}
A.~N{\'e}methi.
\newblock On the {Ozsv{\'a}th-Szab{\'o}} invariant of negative definite plumbed
  $3$-manifolds.
\newblock math.GT/0310083.

\bibitem{OwensStrle}
B.~Owens and S.~Strle.
\newblock Rational homology spheres and four-ball genus.
\newblock math.GT/0308073, 2003.

\bibitem{HolDiskTwo}
P.~S. Ozsv{\'a}th and Z.~Szab{\'o}.
\newblock Holomorphic disks and three-manifold invariants: properties and
  applications.
\newblock To appear in {\em Annals of Math.}. math.SG/0105202.

\bibitem{HolDisk}
P.~S. Ozsv{\'a}th and Z.~Szab{\'o}.
\newblock Holomorphic disks and topological invariants for closed
  three-manifolds.
\newblock To appear in {\em Annals of Math.}. math.SG/0101206.

\bibitem{HolDiskFour}
P.~S. Ozsv{\'a}th and Z.~Szab{\'o}.
\newblock Holomorphic triangles and invariants for smooth four-manifolds.
\newblock math.SG/0110169.

\bibitem{UnknotOne}
P.~S. Ozsv{\'a}th and Z.~Szab{\'o}.
\newblock Knots with unknotting number one and {H}eegaard {F}loer homology.
\newblock math.GT/0401426.

\bibitem{NoteLens}
P.~S. Ozsv{\'a}th and Z.~Szab{\'o}.
\newblock On knot {F}loer homology and lens space surgeries.
\newblock math.GT/0303017.

\bibitem{HolDiskContact}
P.~S. Ozsv{\'a}th and Z.~Szab{\'o}.
\newblock Heegaard {F}loer homologies and contact structures.
\newblock math.SG/0210127, 2002.

\bibitem{HolDiskSymp}
P.~S. Ozsv{\'a}th and Z.~Szab{\'o}.
\newblock Holomorphic disk invariants for symplectic four-manifolds.
\newblock math.SG/0201059, 2002.

\bibitem{HolDiskKnots}
P.~S. Ozsv{\'a}th and Z.~Szab{\'o}.
\newblock Holomorphic disks and knot invariants.
\newblock math.GT/0209056, 2002.

\bibitem{AbsGraded}
P.~S. Ozsv{\'a}th and Z.~Szab{\'o}.
\newblock Absolutely graded {F}loer homologies and intersection forms for
  four-manifolds with boundary.
\newblock {\em Adv. Math.}, 173(2):179--261, 2003.

\bibitem{AltKnots}
P.~S. Ozsv{\'a}th and Z.~Szab{\'o}.
\newblock Heegaard {F}loer homology and alternating knots.
\newblock {\em Geom. Topol.}, 7:225--254 (electronic), 2003.

\bibitem{GenusBounds}
P.~S. Ozsv{\'a}th and Z.~Szab{\'o}.
\newblock Holomorphic disks and genus bounds.
\newblock math.GT/0311496, 2003.

\bibitem{FourBall}
P.~S. Ozsv{\'a}th and Z.~Szab{\'o}.
\newblock Knot {F}loer homology and the four-ball genus.
\newblock {\em Geom. Topol.}, 7:615--639, 2003.

\bibitem{calcKT}
P.~S. Ozsv{\'a}th and Z.~Szab{\'o}.
\newblock Knot {F}loer homology, genus bounds, and mutation.
\newblock math.GT/0303225, 2003.

\bibitem{SomePlumbs}
P.~S. Ozsv{\'a}th and Z.~Szab{\'o}.
\newblock On the {F}loer homology of plumbed three-manifolds.
\newblock {\em Geometry and Topology}, 7:185--224, 2003.

\bibitem{BrDCov}
P.~S. Ozsv{\'a}th and Z.~Szab{\'o}.
\newblock On the {H}eegaard {F}loer homology of branched double-covers.
\newblock math.GT/0309170, 2003.

\bibitem{PlamenevskayaTB}
O.~Plamenevskaya.
\newblock Bounds for {Thurston-Bennequin} number from {Floer} homology.
\newblock math.SG/0311090, 2003.

\bibitem{Rasmussen}
J.~A. Rasmussen.
\newblock Floer homology of surgeries on two-bridge knots.
\newblock {\em Algebr. Geom. Topol.}, 2:757--789 (electronic), 2002.

\bibitem{RasmussenThesis}
J.~A. Rasmussen.
\newblock {\em Floer homology and knot complements}.
\newblock PhD thesis, Harvard University, 2003.
\newblock math.GT/0306378.

\bibitem{RasmussenKhovanov}
J.~A. Rasmussen.
\newblock Khovanov homology and the slice genus.
\newblock Preprint, 2003.

\bibitem{RobbinSalamon}
J.~Robbin and D.~Salamon.
\newblock The {M}aslov index for paths.
\newblock {\em Topology}, 32(4):827--844, 1993.

\bibitem{Rudolph}
L.~Rudolph.
\newblock Quasipositivity as an obstruction to sliceness.
\newblock {\em Bull. Amer. Math. Soc. (N.S.)}, 29(1):51--59, 1993.

\bibitem{Rustamov}
R.~Rustamov.
\newblock Calculation of {H}eegaard {F}loer homology for a class of {B}rieskorn
  spheres.
\newblock math.SG/0312071, 2003.

\bibitem{AFSalamon}
D.~Salamon.
\newblock Lagrangian intersections, {$3$}-manifolds with boundary, and the
  {Atiyah-Floer} conjecture.
\newblock In {\em Proceedings of the International Congress of Mathematicians},
  pages 526--536. Birkh{\"a}user, 1994.

\bibitem{Saveliev}
N.~Saveliev.
\newblock {\em Lectures on the topology of {$3$}-manifolds}.
\newblock de Gruyter Textbook. Walter de Gruyter \& Co., Berlin, 1999.
\newblock An introduction to the Casson invariant.

\bibitem{SeidelDehn}
P.~Seidel.
\newblock The symplectic {F}loer homology of a {D}ehn twist.
\newblock {\em Math. Res. Lett.}, 3(6):829--834, 1996.

\bibitem{SeidelMutation}
P.~Seidel.
\newblock Vanishing cycles and mutation.
\newblock In {\em European Congress of Mathematics, Vol. II (Barcelona, 2000)},
  volume 202 of {\em Progr. Math.}, pages 65--85. Birkh\"auser, Basel, 2001.

\bibitem{SeidelExactSeq}
P.~Seidel.
\newblock A long exact sequence for symplectic {F}loer cohomology.
\newblock {\em Topology}, 42(5):1003--1063, 2003.

\bibitem{Singer}
J.~Singer.
\newblock Three-dimensional manifolds and their {H}eegaard diagrams.
\newblock {\em Trans. Amer. Math. Soc.}, 35(1):88--111, 1933.

\bibitem{TaubesCasson}
C.~H. Taubes.
\newblock Casson's invariant and gauge theory.
\newblock {\em J. Differential Geom.}, 31(2):547--599, 1990.

\bibitem{TaubesSympI}
C.~H. Taubes.
\newblock The {S}eiberg-{W}itten invariants and symplectic forms.
\newblock {\em Math. Research Letters}, 1(6):809--822, 1994.

\bibitem{TaubesSympII}
C.~H. Taubes.
\newblock More constraints on symplectic forms from {S}eiberg-{W}itten
  invariants.
\newblock {\em Math. Research Letters}, 2(1):9--13, 1995.

\bibitem{TaubesBook}
C.~H. Taubes.
\newblock {\em Seiberg {W}itten and {G}romov invariants for symplectic
  {$4$}-manifolds}, volume~2 of {\em First International Press Lecture Series}.
\newblock International Press, Somerville, MA, 2000.
\newblock Edited by R. Wentworth.

\bibitem{TaubesApproach}
C.~H. Taubes.
\newblock Seiberg-{W}itten invariants, self-dual harmonic 2-forms and the
  {H}ofer-{W}ysocki-{Z}ehnder formalism.
\newblock In {\em Surveys in differential geometry}, Surv. Differ. Geom., VII,
  pages 625--672. Int. Press, Somerville, MA, 2000.

\bibitem{Turaev}
V.~Turaev.
\newblock Torsion invariants of {S}pin{$^c$}-structures on $3$-manifolds.
\newblock {\em Math. Research Letters}, 4:679--695, 1997.

\bibitem{Wehrheim}
K.~Wehrheim.
\newblock {\em Anti-self-dual instantonts with {L}agrangian boundary
  conditions}.
\newblock PhD thesis, Swiss Federal Institute of Technology, 2002.

\bibitem{Witten}
E.~Witten.
\newblock Monopoles and four-manifolds.
\newblock {\em Math. Research Letters}, 1:769--796, 1994.

\end{thebibliography}
}

\end{document}